\documentclass[11pt]{article}

\usepackage{multirow, multicol, color, booktabs }
\usepackage{booktabs} % for much better looking tables
\usepackage{array} % for better arrays (eg matrices) in maths
\usepackage{paralist} % very flexible & customisable lists (eg. enumerate/itemize, etc.)
\usepackage{verbatim} % adds environment for commenting out blocks of text & for better verbatim
\usepackage{subfig} % make it possible to include more than one captioned figure/table in a single float
\usepackage{amsmath, amssymb, enumerate}

\newtheorem{theorem}{Theorem}
\newtheorem{corollary}[theorem]{Corollary}
\newtheorem{remark}[theorem]{Remark}
\newtheorem{lemma}[theorem]{Lemma}
\newtheorem{definition}[theorem]{Definition}
\newtheorem{proposition}[theorem]{Proposition}
\newtheorem{example}[theorem]{Example}

\newcommand{\proof}{ {\sc Proof.\quad}}
\newcommand{\pend}{ \hfill $\square$ \\}

\numberwithin{equation}{section}  % Formeln mit f¸hrender Kapitelnummer
\numberwithin{figure}{section}    % Abbildungen mit f¸hrender Kapitelnummer
\numberwithin{table}{section}     % Tabellen mit f¸hrender Kapitelnummer
\numberwithin{theorem}{section}

\newcommand{\of}[1]{\ensuremath{\left( #1 \right)}}

\newcommand{\abs}[1]{\ensuremath{\left| #1 \right|}}
\newcommand{\cb}[1]{\ensuremath{ \left\{ #1 \right\} }}
\newcommand{\sqb}[1]{\ensuremath{ \left[ #1 \right] }}

\newcommand{\bs}{\backslash}

\newcommand{\st}{\,|\;}

\newcommand{\eps}{\ensuremath{\varepsilon}}

\newcommand{\vp}{\ensuremath{\varphi}}

\newcommand{\R}{\mathrm{I\negthinspace R}}
\newcommand{\OLR}{\overline{\mathrm{I\negthinspace R}}}
\newcommand{\N}{\mathrm{I\negthinspace N}}

\renewcommand{\P}{\ensuremath{\mathcal{P}}}

\newcommand{\G}{\ensuremath{\mathcal{G}}}

\newcommand{\C}{\ensuremath{C^-\bs\negthinspace\cb{0}}}

\newcommand{\Min}{{\rm Min}}

\newcommand{\dom}{{\rm dom \,}}

\newcommand{\gr}{{\rm graph \,}}
\newcommand{\cl}{{\rm cl \,}}

\newcommand{\co}{{\rm co \,}}
\newcommand{\cone}{{\rm cone\,}}

\newcommand{\Int}{{\rm int\,}}
\newcommand{\extd}{{\rm extd\,}}

\newcommand{\isum}{{+^{\negmedspace\centerdot\,}}}

\newcommand{\idif}{{-^{\negmedspace\centerdot\,}}}

\newcommand{\lel}{\preccurlyeq}
\newcommand{\leu}{\curlyeqprec}

\usepackage[%backref=page,
colorlinks=true, linkcolor=blue, citecolor=blue, urlcolor=blue, filecolor=blue
]{hyperref}
\definecolor{color0}{gray}{.50}
\definecolor{color1}{rgb}{0,.2,.8}
\definecolor{color2}{rgb}{1,.2,0}
\definecolor{color3}{rgb}{.8,.5,1}

\begin{document}
\title{A Minty variational principle for set optimization}

\author{
Giovanni P. Crespi
\thanks{University of Valle d'Aosta,
Department of Economics and Politics, Loc. Grand Chemin 73-75, 11020 Saint-Christophe, Aosta, Italy.
\href{mailto:g.crespi@univda.it}{g.crespi@univda.it}}
\and
Andreas H. Hamel
\thanks{Free University Bozen--Bolzano,
School of Economics and Management, Universit\"atsplatz - Piazza Universit\'{a} 1, 39031 Bruneck-Brunico
\href{mailto: andreas.hamel@unibz.it}{andreas.hamel@unibz.it}}
\and
Carola Schrage
\thanks{University of Valle d'Aosta,
Department of Economics and Politics, Loc. Grand Chemin 73-75, 11020 Saint-Christophe, Aosta, Italy.
\href{mailto:carolaschrage@gmail.com}{carolaschrage@gmail.com}}
}
\date{{\small \today}}

\maketitle
\begin{abstract}
Extremal problems are studied involving an objective function with values in (order) complete lattices of sets generated by so called set relations. Contrary to the popular paradigm in vector optimization, the solution concept for such problems, introduced by F. Heyde and A. L\"ohne, comprises the attainment of the infimum as well as a minimality property. The main result is a Minty type variational inequality for set optimization problems which provides a sufficient optimality condition under lower semicontinuity assumptions and a necessary condition under appropriate generalized convexity assumptions. The variational inequality is based on a new Dini directional derivative for set-valued functions which is defined in terms of a ``lattice difference quotient": A residual operation in a lattice of sets replaces the inverse addition in linear spaces. Relationships to families of scalar problems are pointed out and used for proofs: The appearance of improper scalarizations poses a major difficulty which is dealt with by extending known scalar results such as Diewert's theorem to improper functions.
\\

{\bf Keywords:} Variational inequalities, Set optimization, Generalized convexity, Dini derivative, residuation.
\end{abstract}

%\input{Intro}
%%%New section
\section{Introduction}

Throughout the paper, let $X$ and $Z$ be two locally convex, topological linear spaces and $C \subseteq Z$ a convex cone with $0 \in C$. Moreover, $\mathcal P(Z)$ denotes the set of all subsets of $Z$ including $\emptyset$. Let a function $f \colon X \to \mathcal P(Z)$ be given. The basic problem is
\[
 \text{minimize} \quad f \quad \text{subject to} \quad x \in X.
\]
Motivated by duality for vector optimization, such set-valued optimization problems have been considered first by Corley  \cite{Corley87JOTA,Corley88JOTA} and Dinh The Luc \cite{Luc89Book}. They gained popularity after the appearance of \cite{KuroiwaTanakaTruong97} and  \cite{Kuroiwa97, Kuroiwa98-1, Kuroiwa98-2} in which so-called set relations are investigated and used to define minimality concepts for sets.

However, the power set $\mathcal P(Z)$ is too large an object and lacks reasonable structure which can be exploited for optimization purposes. On the other hand, additional assumptions imposed to $f$ often imply that the images of $f$ belong to a relatively small subset of $ \mathcal P(Z)$ which carries a richer algebraic and order structure. For example, $C$-convexity of $f$ (see \cite[Definition 1.1]{Borwein77MP}) implies that the set $f(x) + C$ is convex for all $x \in X$. Therefore, appropriate subsets of $\mathcal P(Z)$ are used as image sets of set-valued functions, for example in \cite{Hamel09,HamelSchrage12,Loehne11Book,Schrage10Opt}, and we will follow this approach. The main goal is to define new lower directional derivatives of Dini type for set-valued functions and provide necessary and sufficient conditions in terms of variational inequalities of Minty type to characterize solutions of set-valued minimization problems.

Two questions arise. First, what is understood by a solution of the above problem? Secondly, how can a directional derivative, in particular a difference quotient, be defined if the image set of the function is not a linear space? The answer to the first question is a new solution concept for set-valued optimization problems proposed by F. Heyde and A. L\"ohne \cite{HeydeLoehne11,Loehne11Book}. This concept subsumes classical minimality notions borrowed from vector optimization as well as the infimum/supremum in complete lattices (which are usually not present in vector optimization). The answer to the second is provided by means of residuation operations in (order) complete lattices of sets which replace the inverse addition (the difference) in linear spaces. This approach has been proposed in \cite{HamelSchrage12,HamelSchrage12ArX}. 

Several notions of derivatives for set-valued functions have been introduced, compare e.g. \cite{AlonsoRodriguezMarin2005NonlinAn,AubinFrankowska90,CrespiGinRoc2005,pjov6n1cre,demyanov1986approximation,JahnRau97,KhanRaciti03,Mordukhovich06,Yang98} to mention but a few. Apart from approaches relying on an embedding procedure into a linear space or approaches similar to those in \cite{CrespiGinRoc2005,KhanRaciti03,Yang98}, usually some kind of tangent cone to the graph of $f$ at a point $(x,z) \in X\times Z$ with $z\in f(x)$ is defined to be the graph of the derivative. In this paper, we define a set-valued derivative using increments of function values where the difference is replaced by a residual operation and thus provides a substitute for the the difference quotient in linear spaces. A ``lattice limit" procedure then provides the desired derivative.

It turns out that the lattice concepts are appropriate and sufficient to formulate Minty type variational inequalities which yield the desired characterizations for the new type of solutions. Minty variational inequalities have been introduced in \cite{Minty67} as the problem of finding some $\bar x \in K$ such that
$$
\forall y\in K \colon \langle F\of{y}, \bar x-y\rangle \le 0,
$$
where $F \colon \R^n \to \R^n$, and $K \subseteq \R^n$ is a non empty convex subset. This inequality proved to be useful to study primitive optimization problems when $F$ is some derivative of the objective function $f \colon \R^n\to\R$. The main result in this field is known as Minty variational principle and basically states that the Minty variational inequality provides a sufficient optimality condition for minimizers of $f$ under a lower semicontinuity assumption. The same inequality is also a necessary optimality condition under generalized convexity type assumptions. In \cite{crespiJOTA}, the Minty variational principle has been applied to a non-differentiable scalar optimization problem using lower Dini derivatives. The same approach has been extended to the vector case in \cite{CrespiGinRoc2005}.

The main purpose of this paper is to provide a Minty variational principle for set optimization problems. In the process we also need to deepen the study of lower semicontinuity and generalized convexity. Indeed, it turns out that known results on generalized convexity need to be extended to cover the case of improper functions, which is, to the best of our knowledge, not covered by the existing literature.

The paper is organized as follows. In Section 2, basic notation and results on the ``lattice approach" to set optimization are introduced. The notion of a conlinear space as a natural setting for the image space of classes of set-valued functions is presented in Subsection 2.2. The solution concept for set optimization problems and scalarization techniques are described subsequently. In Section 3, the Dini-type derivative for set-valued functions is introduced, while in Section 4 generalized convexity concepts for possibly improper scalar and set-valued functions are discussed. The main results are presented in Sections 5 and 6 which provide the desired optimality conditions of Minty type for set optimization problems. In Section 7, conclusions are drawn which tie the previous results into a Minty variational principle for set-valued functions.

%%%New section
\section{Functions mapping into complete lattices of sets}

%%%New subsection
\subsection{Some standard notation}

A set $C \subseteq Z$ is called a cone if $z \in C$ and $t>0$ imply $tz \in C$, thus a cone does not necessarily include zero. The conical hull of $C \subseteq Z$ is the set
\[
\cone C = \cb{tz \mid t > 0,\; z \in C}.
\]
The effective domain of an extended real-valued function $\vp \colon X \to \OLR = \R\cup\cb{\pm\infty}$ is the set
$\dom \vp =\cb{x \in X \mid \vp(x) \neq +\infty}$. The lower level sets of such a function are
\[
L_\vp(r)= \cb{x \in X \mid \vp(x) \leq r}, \; r \in \OLR.
\]
This means $\dom \vp=\bigcup\limits_{r\in\R}L_\vp(r)$ and $L_\vp(-\infty)=\bigcap\limits_{r\in\R}L_\vp(r)$. It is well-known that $\vp$ is l.s.c. if, and only if, each lower level set $L_\vp(r)$ with $r \neq +\infty$ is closed.  In this case $L_\vp(-\infty)$ is a closed set.

%%%New subsection
\subsection{The image space}

Order complete lattices of sets which will serve as image spaces for set-valued optimization problems can be generated as follows. The Minkowski (element-wise) addition for non-empty subsets of $Z$ is extended to $\P\of{Z}$ by $\emptyset + A = A + \emptyset = \emptyset$
for $A \in \P\of{Z}$. We shall also write $z + A$ for $\cb{z} + A$ and $z - A$ for $z + (-1)A$ with $-A = \cb{-a \mid a \in A}$. By
\begin{equation}
\label{EqVectorOrder}
z_1 \leq_C z_2 \; \Leftrightarrow \; \cb{z_2} + C \subseteq \cb{z_1} + C  \; \Leftrightarrow \; \cb{z_1} - C \subseteq \cb{z_2} - C
\end{equation}
a preorder (a reflexive and transitive relation) on $Z$ is defined. It is compatible with the linear structure on $Z$, and it is antisymmetric (hence a partial order) if, and only if, $C$ is pointed, that is $C \cap \of{-C} = \cb{0}$. The relation $\leq_C$ can be extended to an order relation on $\P\of{Z}$ in two ways,  defining $A \lel_C B$ $\Leftrightarrow$ $B \subseteq A + C$ and $A \leu_C B$ $\Leftrightarrow$ $A \subseteq A - C$ for $A, B \in \mathcal P(Z)$. These two relations on $\P\of{Z}$ are not antisymmetric in general, and they are different.

One may observe that $A \lel_C B$ if, and only if, $A + C \supseteq B + C$. Thus, on $\cb{A \subseteq Z \mid A = A + C}$ the relation $\lel_C$ is a partial order which coincides with $\supseteq$. If one needs to require closedness and convexity, one is led to
\[
\G\of{Z,C} = \cb{A \subseteq Z \mid A=\cl\co\of{A+C}}
\]
which will be used as an image space for set-valued functions in this note. Here, $\cl A$ and $\co A$ denote the closure and convex hull of a $A \subseteq Z$. Therefore, the relation $\supseteq$ can be understood as an extension of $\leq_C$ from $Z$ to $\G\of{Z,C}$.

For further motivation and more details we refer for example to \cite{Hamel05Habil,Hamel09,Schrage10Opt}. The relation $\lel_C$ is one of the two popularized by Kuroiwa (see e.g.  \cite{Kuroiwa98-1,Kuroiwa98-2,KuroiwaTanakaTruong97}), who originally used them to define solution concepts for optimization problems with a set-valued objective function. The same order relations are applied earlier in other contexts, see e.g. \cite{Brink93AU} and the references in \cite{Hamel05Habil,Jahn04Book}.

A basic observation is as  follows. The pair $\of{\G\of{Z,C}, \supseteq}$ is an order complete, partially ordered set. If $\mathcal A \subseteq \G\of{Z, C}$, then
\begin{equation}
\label{EqInfSupG}
\inf\mathcal A = \cl\co\bigcup\limits_{A \in \mathcal A}A \quad \text{and} \quad \sup\mathcal A = \bigcap\limits_{A\in\mathcal A}A
\end{equation}
as a straightforward check may show. One may also observe $\G\of{Z, C} = \G\of{Z, \cl C}$, so we assume $C = \cl C$ in the following. Therefore, in the remainder of the paper $C \subseteq Z$ is a closed convex cone.

We will also make use of minimal elements with respect to set orders. An element $\bar A \in \mathcal A$ is called a minimal element of $\mathcal A  \subseteq  \G\of{Z,C}$ if
\[
A \in \mathcal A, \; A \supseteq \bar A \quad \Rightarrow \quad A = \bar A.
\]
The set of minimal elements of $\mathcal A$ is denoted by $\Min\of{\mathcal A, \supseteq}$.

Modifying the Minkowski sum and multiplication with nonnegative reals by setting $A \oplus B = \cl\of{A+B}$ and $0 \cdot A = C$ for all $A, B \in \G(Z,C)$ we obtain that $(\G(Z,C),\oplus,\cdot)$ is a (real) conlinear space, i.e. $(\G(Z,C),\oplus,\cdot)$ is a commutative monoid with neutral element $C$ and for all $A, A_1, A_2 \in \G(Z,C)$, $r, s \in \R_+$ it holds $r(A_1\oplus A_2)=rA_1\oplus rA_2$, $r(sA)=(rs)A$ and $1A=A$, $0A=C$, compare \cite[Section 2.1.2]{Hamel05Habil} and \cite{Hamel09}.

Moreover, the order $\supseteq$ on $\G(Z, C)$ is compatible with the algebraic structure of $\of{\G(Z,C),\oplus,\cdot}$, thus $\of{\G(Z,C),\oplus,\cdot, \supseteq}$ is an ordered conlinear space in the sense of \cite[Section 2.1.2]{Hamel05Habil} and \cite{Hamel09}. Finally,
\[
\forall A \in \G(Z,C), \; \forall \mathcal B \subseteq\G(Z,C) \colon A \oplus \inf\mathcal B = \inf(A \oplus \mathcal B)
\]
which provides another link between the algebraic and the order structure: $\of{\G(Z,C),\oplus,\cdot, \supseteq}$ is an inf-residuated conlinear space: compare \cite{GetanMaLeSi,HamelSchrage12, HamelSchrage12ArX}  and the references therein on $\inf$-residuated sets.

The $\inf$-residual of two elements  $A, B \in \G(Z,C)$ is given by
\begin{align*}
A \idif B =\inf\cb{D \in \G(Z,C) \mid B\oplus D\subseteq A} = \cb{z \in Z \mid B + z \subseteq A}.
\end{align*}
For the last equation, compare \cite{HamelSchrage12,HamelSchrage12ArX}. The inf-residual may be seen as a replacement for the inverse addition (the difference) in linear spaces. Indeed, if $A = \cb{z_A} \oplus C = z_A + C$ and $B = \cb{z_B} \oplus C = z_B + C$ then
\[
A \idif B =  \cb{z \in Z \mid  z_B + C + z \subseteq \cb{z_A} \oplus C} = z_A - z_B + C = \cb{z_A - z_B} \oplus C.
\]
A whole calculus for residuals exists, see \cite{HamelSchrage12ArX}. For example, $B\oplus(A\idif B)\subseteq A$ whenever $A, B \in \G(Z,C)$. Compare \cite{HamelSchrage12,HamelSchrage12ArX} and also \cite[Section 4]{penot2011} and the references therein on the use of the residual of two sets.
The concept of residuation, rarely used in (convex) analysis, dates back to Dedekind, \cite[p. 329-330]{Dedekind1872}, \cite[p. 71]{Dedekind1863}, see also
\cite{Birkhoff40, Fuchs66, galatos2007residuated,GetanMaLeSi}.

\begin{example}
\label{ExExtReals}
Let us consider $Z = \R$, $C = \R_+$. Then $\G\of{Z, C} = \cb{[r, +\infty) \mid r \in \R}\cup\cb{\R}\cup\cb{\emptyset}$, and $\G\of{Z, C}$ can be identified (with respect to the algebraic and order structures as introduced above which turn $\G\of{\R, \R_+}$ into an ordered conlinear space and a  complete lattice admitting an inf-residuation) with $\OLR = \R\cup\cb{\pm\infty}$ using the 'inf-addition' $\isum$ (see \cite{HamelSchrage12,RockafellarWets98}) and the inf-residuation on $\OLR$ as given by
\[
r \isum s = \inf\cb{a+b \mid a, b \in \R, \, r \leq a, \, s \leq b} \quad \text{and} \quad
r \idif s  = \inf\cb{t \in \R \mid r \leq s \isum t}
\]
for all $r,s\in\OLR$, compare \cite{HamelSchrage12,HamelSchrage12ArX} for further details.
\end{example}

Simple examples show that the inf-residual of two sets $A, B \in \G\of{Z, C}$ can be empty in many (interesting) cases. One may realize what is going on taking $Z=\R^2$, $C = \R^2_+ = A$ and $B = \cb{z \in \R^2 \mid 2z_1 + z_2 \geq 0, \; \frac{1}{2}z_1 + z_2 \geq 0}$. Therefore, we introduce another such operation.

Let $Z^*$ be the topological dual of $Z$. The (negative) dual cone of $C$ is the set
\[
C^- = \cb{z^* \in Z^* \mid \forall c \in C \colon z^*(c) \leq 0}.
\]
We assume $\C \neq \emptyset$. Take $z^*\in \C$ and define
\[
H\of{z^*} = \cb{z \in Z \mid z^*\of{z} \leq 0}
\]
which is the homogeneous closed half space with normal $z^*$. The $z^*$-residual of $A, B \in \G(Z,C)$ is
\[
A-_{z^*} B = \of{A \oplus H\of{z^*}} \idif B \\
	= \cb{z \in Z \mid B + z \subseteq A \oplus H\of{z^*}},
\]

Of course, $A-_{z^*} B$ coincides with the inf-residual of $A \oplus H\of{z^*}$ and $B \oplus H\of{z^*}$ as elements of $\G\of{Z, H\of{z^*}}$.

%%%New subsection
\subsection{$\G\of{Z, C}$-valued functions and the solution concept}

Let $f \colon X \to \G(Z,C)$ be a function. The graph and the effective domain of $f$ are the sets
\[
\gr f = \cb{(x,z) \in X \times Z \mid z \in f(x)} \quad \text{and} \quad \dom f = \cb{x \in X \mid f(x) \neq \emptyset},
\]
respectively. The function $f$ is called convex if $\gr f$ is convex, and it is called positively homogeneous if $\gr f$ is a cone. The set $f\sqb{M} = \cb{f\of{x} \mid x \in M}$ is the image of $M \subseteq X$ under $f$. In particular, $f\sqb{X}$ is the image of $X$ under $f$.

What shall we understand by a solution of a set-valued optimization problem? The traditional idea is to look for points $\of{\bar x, \bar z} \in \gr f$ such that $\bar z$ is a minimal point of $\bigcup_{x \in X} f\of{x}$ according to $\leq_C$. However, this is not very satisfactory in many cases (see, for example, \cite[p. 210]{JahnTruong11}), and therefore, the so-called set relation approach has been proposed (\cite{KuroiwaTanakaTruong97} and several papers by D. Kuroiwa, among them \cite{Kuroiwa98-1,Kuroiwa98-2}) which consists of looking for $\bar x \in \dom f$ such that the value $f\of{\bar x}$ is minimal with respect to a set relation in the set $f\sqb{X} = \cb{f(x) \mid x \in X}$. The reader may compare \cite{HernandezRodriguezMarin07} and \cite{RodriguezMarinSama09} for this approach.

The paper \cite{HeydeLoehne11} (see also \cite{Loehne11Book}) put forth a new idea which in some way synthesizes the two previous approaches. Here is the basic definition adapted to our framework.

\begin{definition}
\label{DefSetSolution} (a) A set $M \subseteq X$ is called an infimizer of the function $f \colon X \to \G\of{Z, C}$ if
\[
\inf_{m \in M}f(m) = \inf_{x \in X}f(x).
\]

(b) A point $\bar x \in X$ is called a minimizer of $f \colon X \to \G\of{Z, C}$ if $f\of{\bar x}$ is a minimal element of $f\sqb{X}$, i.e.
\[
x \in X, \; f\of{x} \supseteq f\of{\bar x} \quad \Rightarrow  f\of{x} = f\of{\bar x}.
\]

(c) A set $M \subseteq X$ is called a solution of the problem
\begin{align}\label{eq:P}
\tag{P} \text{minimize} \quad f \colon X \to \G\of{Z, C} \quad \text{subject to} \quad x \in X,
\end{align}
if $M$ is an infimizer of $f$ consisting only of minimizers.
\end{definition}

Note that the condition in (a) is equivalent to
\[
\cl\co\bigcup\limits_{m \in M} f\of{m} = \cl\co\bigcup\limits_{x \in X} f\of{x},
\]
and this condition means that the infimum of $f$ is attained in $M$. The condition in (b) just is $f(\bar x) \in \Min\of{f\sqb{X}, \supseteq}$. In the set-valued framework, or as a matter of fact already in the (multi-dimensional) vector-valued framework, infimum attainment and minimality are no longer equivalent (as in the scalar case), but they should, of course, both be part of a meaningful solution concept.
The following examples discuss a few difficulties one encounters when looking at
``vector solutions'' rather then ``set solutions''.

\begin{example}
\label{ExMinVsInf}
(a) Consider $\OLR$ with the usual relation  $\leq$ and $\vp \colon X \to\OLR$. Then $\vp\sqb{X}$ possesses minimal elements, if and only if $\inf\vp\sqb{X} \in \vp\sqb{X}$, and $\inf \vp\sqb{M}=\inf\vp\sqb{ X}$ is a solution of \eqref{eq:P}, iff $\vp(m)=\inf\cb{\vp(x) \mid x \in X}$ for all $m \in M$.

(b) Let $S = \R^2\cup\cb{\pm\infty}$ be equipped with the order $\leq_C$ generated by the convex cone $C = \R^2_+$ with the obvious extension to $\pm\infty$. Define a function $F \colon X = \R^2 \to S$ by
\[
F\of{x} = \left\{
	\begin{array}{ccl}
		x & : & \text{ if } 1 - x_1 \leq x_2, \; 0 \leq x_1 \leq 1; \\
		+\infty & :&  \text{otherwise.}
	\end{array}
	\right.
\]
The set of minimal elements of $F\sqb{X}$ is the set $M=\cb{F(x) \in \R^2 \mid x_2 = 1 - x_1, \; 0 \leq x_1 \leq 1}$ whereas $\inf\cb{F\of{x} \mid x \in X} = 0 \in \R^2$ does not belong to the range of $F$. 

(c) Let $\of{S, \leq}$ be as in (b) and consider
\[
F\of{x} = \left\{
	\begin{array}{ccl}
		x & : & \text{ if }\max\cb{1 - 2x_1, \frac{1}{2} - \frac{1}{2}x_1} \leq x_2, \; x_1 \in \R; \\
		+\infty & : & \text{otherwise.}
	\end{array}
	\right.
\]
The set of minimal elements of $F\sqb{X}$ is the set  $M=\cb{F(x)\in \R^2 \mid \max\cb{1 - 2x_1, \frac{1}{2} - \frac{1}{2}x_1} = x_2}$ whereas $\inf\cb{F\of{x} \mid x \in X} = -\infty$.
\end{example}

As a remedy for the difficulty in the previous example, a vector-valued problem is embedded into a $\G\of{Z, C}$-valued one. Using the device introduced in Definition \ref{DefSetvaluedExtension} below one may see that (subsets of) the set of minimal points of $F\sqb{X}$ with respect to $\leq_C$ indeed provides solutions of the set-valued problem -- in the sense of Definition \ref{DefSetSolution}.

\begin{definition}
\label{DefSetvaluedExtension}
Let $F \colon X \to Z\cup\cb{\pm\infty}$ be an extended vector-valued function. Its $\G(Z,C)$-valued extension $f \colon X \to \G(Z,C)$ is defined by
\[
f\of{x} =
	\left\{
	\begin{array}{ccl}
	Z & : & F\of{x} = -\infty \\
	F\of{x} + C & : & F\of{x} \in Z \\
	\emptyset & : & F\of{x} = +\infty.
	\end{array}
	\right.
\]
\end{definition}

If $f$ is such a $\G(Z,C)$-valued extension, then $F(x_1) \leq_C F(x_2)$ if, and only if, $f(x_1) \supseteq f(x_2)$ for $x_1, x_2 \in X$.
Hence $f(x) \in \Min\of{f\sqb{X},\supseteq}$ if, and only if, $F(x) \in \Min\of{F[X], \leq_C}$. Moreover, $M$ is a solution of \eqref{eq:P} if, and only if, $F[M] \subseteq \Min\of{F[X], \leq_C}$ and $F\sqb{X}\subseteq \cl\co\of{F\sqb{M}+C}$. The last inclusion can be understood as a weakened version of the so called domination property.

%%%New subsection
\subsection{Scalarizations}

We will associate to $f \colon X \to \G(Z, C)$ a family of extended real-valued functions which we call scalarizations of $f$. For $z^* \in \C$ and $r \in \OLR$, the set
\[
L_{z^*}(-r) = \cb{z \in Z \mid r\leq -z^*(z)}
\]
is a closed half space with normal direction $z^*$ if $r\in\R$, while $L_{z^*}(+\infty)=\emptyset$ and $L_{z^*}(-\infty)=Z$. We have $L_{z^*}(0) = H\of{z^*}$.

\begin{definition}
Let $f \colon X \to \G(Z,C)$ and $z^*\in\C$. The scalarization of $f$ with respect to $z^*$  is the function $\vp_{f,z^*} \colon X \to \OLR$ defined by
\[
\vp_{f,z^*}(x) = \inf\cb{-z^*(z) \mid z \in f(x)}.
\]
\end{definition}

Of course, $-\vp_{f,z^*}(x) = \sup_{z \in f\of{x}} z^*(z)$ is the value of the support function of the set $f\of{x}$ at $z^*$. Since the values of $f$ are closed convex, they are the intersections of all closed half spaces including them; such half spaces can only be generated by elements of $\C$ since $f$ maps into $\G\of{Z, C}$. Hence
\begin{align} \label{EqSetification}
f(x)  & = \bigcap\limits_{z^* \in \C} \cb{z \in Z \mid \vp_{f,z^*}(x) \leq -z^*\of{z}}  =  \bigcap\limits_{z^* \in \C} f(x)\oplus H\of{z^*},
\end{align}
and we have $f(x) \oplus H\of{z^*} = L_{z^*}(-\vp_{f,z^*}(x))$ as well as $\dom f = \dom \vp_{f,z^*}$ for all $z^*\in\C$.

\begin{example}
\label{ExSetvaluedExtension}
Let $F \colon X \to Z\cup\cb{\pm\infty}$ be an extended vector-valued function and $f \colon X \to \G(Z,C)$ its $\G(Z,C)$-valued extension.
Then, $\dom f = \cb{x \in X \mid F\of{x} \neq +\infty}$, and for each $z^*\in \C$, $\vp_{f, z^*}(x) = -\of{z^* \circ F}(x) = -z^*\of{F\of{x}}$ whenever $F\of{x} \in Z$, $\vp_{f,z^*}(x)=-\infty$ whenver $F\of{x} = -\infty$ and $\vp_{f,z^*}(x)=+\infty$ whenver $F\of{x} = +\infty$.
\end{example}

\begin{remark}
\label{RemSolSetValuedScalarized}
If $f \colon X \to \G\of{Z,C}$, then a set $M\subseteq X$ is a solution to \eqref{eq:P} if, and only if,

(a) $\forall x \in X$: $f\of{x} \subseteq \cl\co\bigcup\limits_{m\in M}f(m)$,

(b) for $m\in M$ and $x \in X$, either $\vp_{f,z^*}(m) = \vp_{f,z^*}(x)$
for all $z^*\in\C$, or there exists $z^*_0\in\C$ such that $\vp_{f,z_0^*}(m) < \vp_{f,z_0^*}(x)$.

\medskip This follows by a standard separation argument since the requirements for a solution as given in Definition \ref{DefSetSolution} only leave two possibilities for $f\of{x}, f\of{m} \in \G\of{Z,C}$: either $f\of{x} = f\of{m}$, or $f\of{x} \not\supseteq f\of{m}$.
\end{remark}

%%%New subsection
\subsection{Continuity notions for set-valued functions}

By $\mathcal U_X(0)$ and $\mathcal U_Z(0)$ we denote a neighborhood base at the origin of $X$ and $Z$, respectively. If necessary, we assume that such a neighborhood base consists of absorbing, balanced and convex sets which is always possible in locally convex spaces.

\begin{definition} \label{DefLsc} (a) A function $f \colon X \to \G\of{Z,C}$ is said to be lattice-lower semicontinuous (lattice-l.s.c. for short) at $x \in X$ if
\[
\forall x \in X \colon f(x)\leq \liminf\limits_{y \to x}f(y)  = \bigcap_{U \in \mathcal U_X\of{0}} \cl\co \bigcup_{u \in U} f\of{x + u}.
\]
It is called lattice-lower semicontinuous if it is lattice-lower semicontinuous at each $x \in X$.

(b) A function $f \colon X \to \G\of{Z,C}$ is said to be uniformly (Hausdorff) lower semicontinuous (uniformly-l.s.c. for short) if for all  $V \in \mathcal U_Z(0)$ there exists a $U \in \mathcal U_X(0)$ such that for all $x \in X$
\begin{equation}
\label{eq:Uniform_lsc}
\forall u \in U \colon f(x+u) \subseteq f(x)+V.
\end{equation}

(c) A function $f \colon X \to \G\of{Z,C}$ is called $D^*$-lower semicontinuous if $\vp_{f,z^*}$ is an extended real-valued lower semi-continuous function for all $z^* \in D^*$ where $D^* \subseteq \C$, and ``$\cb{z^*}$-lower semicontinuous'' will be abbreviated to ``$z^*$-lower semicontinuous''.
\end{definition}

Recall that a scalar function $\vp \colon X \to \OLR$ is l.s.c. if, and only if, $\vp(x) \leq \liminf\limits_{y \to x}\vp(y)$ for all $x\in X$.

\begin{proposition}\label{prop:Unif_lsc}
(a) If $f \colon X \to \G(Z,C)$ is uniformly l.s.c., then each scalarization $\vp_{f,z^*} \colon X \to \OLR$ with $z^* \in \C$ is uniformly l.s.c., i.e. for all $\varepsilon>0$ there exists a $U\in\mathcal U_X(0)$ such that for all $x \in X$
\begin{align}\label{eq:uniform_lsc_scalar}
 \forall u \in U \colon \vp_{f,z^*}(x)\leq \vp_{f,z^*}(x+u) + \varepsilon.
\end{align}

(b) If each scalarization $\vp_{f,z^*} \colon X \to \OLR$ of $f \colon X \to \G\of{Z,C}$ with $z^* \in \C$ is l.s.c., then $f$ is lattice-l.s.c.
\end{proposition}
\proof
(a) If $V = \cb{z \in Z \mid z^*\of{z} \leq \eps}$ then, by definition, there exists a $U \in \mathcal U_X(0)$ such that $f(x+u) \subseteq f(x) \oplus V$, hence $\vp_{f,z^*}(x)\leq \vp_{f,z^*}(x+u)+\varepsilon$ for all $u\in U$ and all  $x\in X$.

(b) Lower semicontinuity of $\vp_{f,z^*}$ is
\[
\forall x \in X \colon \vp_{f,z^*}(x) \leq \liminf_{y \to x} \vp_{f,z^*}(x)  = \sup_{U \in \mathcal U_X(0)}\inf_{u \in U}\vp_{f,z^*}(x+u).
\]
Since the support function of a set coincides with the support function of the closed convex hull of the same set we obtain
\[
\inf_{u\in U}\vp_{f,z^*}(x+u) = \inf_{u\in U}\inf_{z\in f(x+u)}-z^*(z) = \inf\cb{-z^*(z) \mid z \in \cl\co\bigcup\limits_{u\in U}f(x+u)}.
\]
On the other hand, 
\[
\eta \in \bigcap\limits_{U\in\mathcal U_X(0)}\cl\co\bigcup\limits_{u\in U}f(x+u)
\]
implies
\[
\forall U\in\mathcal U_X(0) \colon \inf\cb{-z^*(z) \mid z \in \cl\co\bigcup\limits_{u\in U}f(x+u)} \leq -z^*(\eta),
\]
hence
\[
\sup_{U\in\mathcal U_X(0)}\inf\cb{-z^*(z) \mid z \in \cl\co\bigcup\limits_{u\in U}f(x+u)} \leq -z^*(\eta),
\]
and finally
\[
\sup_{U\in\mathcal U_X(0)}\inf\cb{-z^*(z) \mid z \in \cl\co\bigcup\limits_{u\in U}f(x+u)} \leq 
	\inf\cb{-z^*(z) \mid \bigcap\limits_{U\in\mathcal U_X(0)}\cl\co\bigcup\limits_{u\in U}f(x+u)}.
\]
Thus, 
\begin{multline*}
f(x) \stackrel{\eqref{EqSetification}}{=} \bigcap\limits_{z^* \in \C} \cb{z \in Z \mid \vp_{f,z^*}(x) \leq -z^*\of{z}} \\
\supseteq \bigcap\limits_{z^*\in\C} \cb{z\in Z\st \inf\cb{-z^*(z)\st z\in \bigcap\limits_{U\in\mathcal U_X(0)}\cl\co\bigcup\limits_{u\in U}f(x+u)}\leq -z^*(z)} = \liminf\limits_{y \to x} f(y)
\end{multline*}
where the last equation also is \eqref{EqSetification} applies to $\liminf\limits_{y \to x} f(y)$ instead of $f(x)$. \pend

The property defined in (b) of Definition \ref{DefLsc} is a uniform version of a continuity notion called Hausdorff upper continuity in \cite{Goepfert03}. We refer to it as ``uniformly-l.s.c.'' to avoid confusion since it implies lattice-lower semicontinuity for a set-valued function and also lower semicontinuity of its scalarizations as shown in (b) of Proposition \ref{prop:Unif_lsc}.

If $f \colon X \to \G(Z,C)$ is the set-valued extension of a function $F \colon X \to Z\cup\cb{\pm\infty}$ as introduced in Definition \ref{DefSetvaluedExtension}, then $f$ is $\C$-l.s.c. if $F$ is (Hausdorff) l.s.c.  in the sense that for all $x \in X$ and for all  $V \in \mathcal U_Z(0)$ there exists a $U \in \mathcal U_X(0)$ such that
\[
\forall u \in U \colon f\of{x+u} \in f\of{x} + V.
\]

For a more detailed comparison among different continuity notions of set-valued functions we refer to \cite{HeydeSchrage13} and the references therein.

%%%New section

\section{The lower Dini directional derivative for set-valued functions}

Here is the definition of a new Dini-type derivative for set-valued functions.
 
\begin{definition}
\label{DefLowerDiniDer}
Let $f \colon X \to \G(Z,C)$ be a function, $x, u \in X$ and $z^*\in \C$. The lower Dini directional derivative of $f$ at $x$ in the direction $u$ with respect to $z^*$ is
\begin{align*}
f^\downarrow_{z^*}(x,u) = \liminf\limits_{t \downarrow 0}\frac{1}{t}\of{f(x+tu) -_{z^*} f(x)}
					 = \bigcap\limits_{s > 0}\cl\bigcup\limits_{t \in \of{0, s}}\frac{1}{t}\of{f(x+tu)-_{z^*}f(x)}.
\end{align*}
\end{definition}

Note that we can drop the convex hull involved in the infimum in $\G(Z,C)$ since the union of closed half spaces with the same normal automatically is convex. For scalar functions we adapt the standard definition of the lower Dini directional derivative to our setting.

\begin{definition}\label{scalarDini}
Let $\vp \colon X \to \OLR$ be a scalar function, $x,u\in X$. The lower 
Dini directional derivative of $\vp$ at $x$ in direction $u$ is
\begin{align*}
\vp^\downarrow(x,u)&=\liminf\limits_{t\downarrow 0}\frac{1}{t}\of{\vp(x+tu)\idif \vp(x)}.
\end{align*}
\end{definition}

With Definition \ref{scalarDini}, we do not restrict to $x\in \dom \vp$, nor we do demand $\vp$ to be a proper function. To this extent, the difference operator is replaced by $\idif$, the residual operator. 

\begin{example}\label{ex:scalar_fct_SV_derivative}
Let $\vp \colon X \to \OLR$ be a extended real-valued function and $f$ its $\G\of{\R, \R_+}$-valued extension (see Definition \ref{DefSetvaluedExtension}). The dual cone of $C=\R_+$ is $-\R_+$ and 
\[
\forall z^* \in \C \colon f^\downarrow_{z^*}(x,u) = f_{-1}^\downarrow(x,u).
\]
Moreover, $f_{-1}^\downarrow\of{x, \cdot}$ is the $\G\of{\R, \R_+}$-valued extension of $\vp^\downarrow\of{x,\cdot}$, that is
\[
f_{-1}^\downarrow\of{x, u}  = 
	\left\{
	\begin{array}{lcl}
	\R & : & \vp^\downarrow\of{x,u} = -\infty \\
	\vp^\downarrow\of{x,u} + \R_+ & : & \vp^\downarrow\of{x,u} \in \R \\
	\emptyset & : & \vp^\downarrow\of{x,u} = +\infty
	\end{array}
	\right.
\]
This can be shown by observing
\begin{align*}
f\of{x + tu} -_{-1} f\of{x} = 
\left\{
	\begin{array}{lcl}
	\R & \Leftrightarrow &  \vp\of{x + tu} \idif \vp\of{x} = -\infty \\ 
	\vp\of{x + tu} - \vp\of{x} + \R_+  & \Leftrightarrow &  \vp\of{x + tu}, \vp\of{x} \in \R \\ 
	\emptyset & \Leftrightarrow & \vp\of{x + tu} \idif \vp\of{x} = +\infty,
	\end{array}
\right.
\end{align*}
and these cases are mutually exclusive and exhausting.
\end{example}

The following proposition collects some elementary properties of Dini derivatives for future reference.

\begin{proposition}
\label{PropDiniElemProps}
(a) Both $u \mapsto f^\downarrow_{z^*}(x,u)$ and $u \mapsto \vp^\downarrow(x,u)$ are positively homogeneous, i.e.
\[
\forall r > 0 \colon f^\downarrow_{z^*}(x,ru) = rf^\downarrow_{z^*}(x,u)
\]
and parallel for $\vp^\downarrow$.

(b) For all $x \in X$, for all $u \in X$,
\begin{align}
\label{EqDirDerScalarized}
f^\downarrow_{z^*}(x,u) & = \cb{z \in Z \mid \vp_{f,z^*}^\downarrow(x,u) \leq -z^*\of{z}} \\
\vp_{f,z^*}^\downarrow(x,u) & = \vp_{f^\downarrow_{z^*}\of{x, \cdot}, z^*}\of{u}.
\end{align}

(c) If $x \notin \dom f = \dom \vp_{f,z^*}$, then $\vp_{f,z^*}^\downarrow(x,u)=-\infty$  and $f^\downarrow_{z^*}(x, u)=Z$ for all $u\in X$.
\end{proposition}

{\sc Proof.} (a) and (c) are immediate. For (b),  observe  that for any $x,u\in X$, $t>0$ and any $z^*\in\C$,
\[
\frac{1}{t}\of{f(x+tu)-_{z^*} f(x)}=\frac{1}{t}\cb{z\in Z\st \vp_{f,z^*}(x+tu)\idif \vp_{f,z^*}(x)\leq -z^*(z)  }
\]
and
\[
\bigcap\limits_{t_0>0}\cl\bigcup\limits_{0<t<t_0}\frac{1}{t}\of{f(x+tu)-_{z^*} f(x)}
=\cb{z\in Z\st \liminf\limits_{t\downarrow 0}\frac{1}{t}\of{\vp(x+tu)\idif \vp(x)}\leq -z^*(z)  },
\]
compare also \cite{HamelSchrage12,Schrage10Opt}.
\pend 

\begin{remark}
\label{RemVectorVsSetDer} This is not the first attempt to introduce a Dini derivative for set-valued functions. In \cite{CrespiGinRoc2005,pjov6n1cre}, for instance, the lower Dini directional derivative of $f$ at $(x,z)$ with $z\in f(x)$ was defined as
\begin{align*}
f'(x,z;u) &= \bigcap\limits_{s > 0}\cl\bigcup\limits_{t\in\of{0,s}}\frac{1}{t}\of{f(x+tu)+\cb{-z}}.
\end{align*}
Since
\begin{multline*}
\forall z \in f\of{x} \colon f(x+tu) -_{z^*} f(x) = \cb{y \in Z \mid f(x) + y \subseteq f(x + tu) \oplus H(z^*)}  \\
	\subseteq \cb{y \in Z \mid z + y \subseteq f(x + tu) \oplus H(z^*)} = \sqb{f(x + tu) + \cb{-z}} \oplus H(z^*)
\end{multline*}
we have $f^\downarrow_{z^*}(x,u) \subseteq f'(x,z;u) \oplus H(z^*) = L_{z^*}\of{\sup\cb{z^*(y) \mid y \in f'(x,z;u)}}$. On the other hand, if there is $z \in f(x)$ such that 
$-z^*(z)=\vp_{f,z^*}(x)$ then one can replace $f(x)$ by $z + H(z^*)$ in the above formula and obtains the converse inclusion, thus  $f^\downarrow_{z^*}(x,u) = f'(x,z;u) \oplus H(z^*)$ in this case.

This means that the $z^*$-Dini derivative is a little more precise than the previous concept which is taken ``at points of the graph''. If an assumption about the existence of support points of $f(x)$ is satisfied then the two concepts coincide ``half space-wise'' at those support points.

In particular, if $f(x)\subseteq z+C$ holds true for $z\in f(x)$, then
\begin{align*}
f'(x,z;u) = \bigcap\limits_{z^*\in\C}f^\downarrow_{z^*}(x,u).
\end{align*}
This shows that for vector-valued functions one can take intersections of the half space-valued $z^*$-Dini derivatives.
\end{remark}

Another idea is to use the residual operation in $\G(Z,C)$ instead of its $z^*$-variant in $\G(Z, H(z^*))$, compare \cite{CrespiSchrage2013SOmF}. The corresponding lower Dini directional derivative of a function $f \colon X \to \G(Z,C)$ is defined by 
\begin{align*}
f^\downarrow(x,u)=\liminf\limits_{t\downarrow 0}\frac{1}{t}\of{f(x+tu)\idif f(x)}
\subseteq  \bigcap\limits_{z^*\in\C}f^\downarrow_{z^*}(x,u).
\end{align*}
The following example shows that this derivative quickly becomes ``non-finite'' in the sense that it assumes the value $\emptyset$ even if the lower Dini derivative with respect to $z^*$ is non-empty for each $z^* \in \C$.

It will become clear in Section \ref{sec:NondomElem} that Definition \ref{DefLowerDiniDer} provides a good enough concept for Minty type variational inequalities.

\begin{example}\label{ex:Strict_domination_IntC_empty} Let $X=\R$ and $Z=\R^2$ with the ordering cone $C=\cl\cone\cb{(0,1)^T}$ and $f \colon X \to \G(Z,C)$ defined by 
\[
f(x)= \left\{
	\begin{array}{lcl}
		\sqb{x^2-1, 1-x^2}\times \R_+  &:&  x\in\sqb{0,1} \\
		\emptyset  &:& \text{otherwise}
	\end{array}
\right.
\]
Fix $x=0$ and $u=1$. Then for all $z^*\in\C$, $t \in \sqb{0,1}$ 
\begin{align*}
\vp_{f,z^*}(x+tu)=\of{t^2-1}\abs{z^*_1} \quad \text{and} \quad \frac{1}{t}\of{\vp_{f,z^*}(x+tu)\idif \vp_{f,z^*}(x)}=t\abs{z^*_1}.
\end{align*}
Thus $\vp^\downarrow_{f,z^*}(x,u)=0$ and $f_{z^*}^\downarrow(x,u)=H(z^*)$.

Let $z=(k,l)^T\in f(x)$ with $k\in\sqb{-1,1}$ and $l= 0$, then

\begin{align*}
f'(x,z;u) &= 
	\left\{
	\begin{array}{lcl}
\cb{(y_1,y_2)^T\in\R^2\st y_2\geq 0}&:&\text{ if } k\neq\pm 1;\\
\cb{(y_1,y_2)^T\in\R^2\st y_1,y_2\geq 0}&:&\text{ if } k=-1;\\
\cb{(y_1,y_2)^T\in\R^2\st -y_1,y_2\geq 0}&:&\text{ if } k=1.\\
	\end{array}
	\right.
\end{align*}
Let $z=(k,l)^T\in f(x)$ with $k\in\sqb{-1,1}$ and $l> 0$, then
\begin{align*}
f'(x,z;u) &= 
	\left\{
	\begin{array}{lcl}
\R^2 &:&\text{ if } k\neq\pm 1;\\
\cb{(y_1,y_2)^T\in\R^2\st y_1\geq 0}&:&\text{ if } k=-1;\\
\cb{(y_1,y_2)^T\in\R^2\st -y_1\geq 0}&:&\text{ if } k=1.\\
	\end{array}
	\right.
\end{align*}

For all $t\in\of{0,1}$, it holds
\begin{multline*}
\frac{1}{t}\of{f(x+tu)\idif f(x)}
=\frac{1}{t}\cb{ z\in Z\st f(x)+z\subseteq f(x+tu)}\\
=\frac{1}{t}\cb{ z\in Z\st \sqb{(-1+z_1,0)^T,(1+z_1,0)^T}\subseteq \sqb{(t^2-1,0)^T,(1-t^2,0)^T}}
=\emptyset,
\end{multline*}
hence 
\[
f^\downarrow(x,u) =\emptyset
\subsetneq \bigcap\limits_{z^*\in\C}L_{z^*}(-\vp^\downarrow_{f,z^*}(x,u))=\bigcap\limits_{z^*\in\C}f_{z^*}^\downarrow(x,u)=C.
\]
\end{example}

%%%New section
\section{Generalized convexity}

Generalized convexity and generalized monotonicity arise almost naturally when dealing with a Minty variational principle (see e.g. \cite{crespiJOTA}).
In the following, we need the following concept.

A set $D \subseteq X$ is said to be star-shaped at $\bar x \in D$ if
\[
\forall x \in D, \; \forall t \in [0,1] \colon t\bar x + (1-t)x \in D.
\]

The results on extended real-valued functions $\vp:X\to\OLR$ presented in the following resemble known results on proper functions, as given e.g. in  \cite{CambiniMartein2009} and even the proofs are in the same line. However, to the best of our knowledge none of the properties or even definitions below has been stated for improper functions, thus proofs are given here for the sake of completeness.

%%%New subsection
\subsection{Extension to the extended real-valued case} %of scalar results

Let $\vp \colon X \to \OLR$ be an extended real-valued function. The function $\vp_{a, b} \colon \R \to \OLR$ is defined by
\[
\vp_{a,b}\of{t} =
	\left\{
	\begin{array}{ccc}
	\vp(a+t(b-a)) & : & t\in\sqb{0,1} \\
	+\infty & : & \text{otherwise}
	\end{array}
	\right.
\]

In the following, we will say that the function $\vp$ is radially l.s.c. at $a$ if the function $\vp_{a, b}$ is l.s.c. for all $b \in X$, and similar for other properties. The following result is Diewert's Mean Value Theorem \cite{Diewert81}.

\begin{proposition}
\label{PropDiewert}
Let $\vp \colon X \to \OLR$ and $a, b\in X$ be such that $\vp_{a,b} \colon \sqb{0,1} \to \R$ is lower semicontinuous (and real-valued). Then, there exist $0 \leq t < 1$ and $0 < s \leq 1$ such that
\begin{align*}%\label{eq:Diewert_MVT}
 \vp(b) -  \vp(a) &\leq ( \vp_{a,b})^\downarrow(t,1) \; \text{and} \\
 \vp(a) -  \vp(b) &\leq ( \vp_{a,b})^\downarrow(s,-1).
\end{align*}
\end{proposition}

Note that for all $0 \leq t < 1$ and $0 < s\leq 1$ the following equations are satisfied
\begin{align*}%\label{eq:Diewert_MVT}
( \vp_{a,b})^\downarrow(t,1)   &=\vp^\downarrow(a+t(b-a),b-a), \\
( \vp_{a,b})^\downarrow(s,-1) &=\vp^\downarrow(a+s(b-a),a-b).
\end{align*}

By a careful case study, we can extend this classical result to the case when $\vp_{a,b} \colon \sqb{0,1}\to \OLR$ is extended real-valued and not necessarily proper. Then, the difference has to be replaced by the inf-residual in $\OLR$.

\begin{theorem}\label{thm:Diewert}
Let $\vp \colon X \to \OLR$ and $a, b \in X$ be given such that $a \neq b$ and $\vp_{a,b} \colon \R \to \OLR$ is lower semicontinuous. Then:

(a) If either $\vp(a) = +\infty$, or $\cb{a,b} \subseteq \dom \vp$, then there exists $0\leq t< 1$ such that
\[
  \vp(b)\idif \vp(a)\leq \of{\vp_{a,b}}^\downarrow(t,1).
\] 

(b) If either $\vp(b)=+\infty$, or $\cb{a,b} \subseteq \dom \vp$, then there exists $0< s \leq 1$ such that
\[
\vp(a) \idif \vp(b) \leq \of{\vp_{a,b}}^\downarrow(s,-1).
\]
\end{theorem}

{\sc Proof.} (a) The proof of the first inequality is given via a case study. If $\vp(a)=+\infty$ or $\vp(b)=-\infty$, then 
\[
\vp(b)\idif\vp(a)=\inf\cb{r\in \R\st \vp(b)\leq \vp(a)+r}=-\infty,
\]
so the first inequality is trivially satisfied. 

Next, assume $\cb{a,b}\subseteq\dom\vp$ and $\vp(b)\neq-\infty$. If $\vp_{a,b}(t)=-\infty$ for some $0\leq t<1$, then by lower semicontinuity $\vp_{a,b}(t_0)=-\infty$, setting
\[
  t_0=\sup\cb{t\in\cb{0,1}\st \vp_{a,b}(t)=-\infty}
\]
and by assumption $t_0<1$. Hence $\of{\vp_{a,b}}^\downarrow(t_0,1)=+\infty$, satisfying the first inequality.

Finally, let $\cb{a,b}\subseteq\dom\vp$ and $\vp(b)\neq-\infty$ be assumed and $\vp_{a,b}(t)=+\infty$ for some $0<t<1$ and set
\[
  t_0=\inf\cb{t\in\of{0,1}\st \vp_{a,b}(t)=+\infty}.
\]
If $t_0=0$, then we are finished, as in this case $\of{\vp_{a,b}}^\downarrow(0,1)=+\infty$ is true, hence assume $0<t_0$. In this case, $\sqb{0,t}\subseteq\dom\vp_{a,b}$ is true for all $t\in\of{0,t_0}$, and the above result combined with Proposition \ref{PropDiewert} applied to $b=a+t(b-a)$ gives that for all $0 < t < t_0$ there exists a $0\leq \bar t<1$ such that
\[
  \vp(a+t(b-a))\leq \vp(a)\isum \of{\vp_{a,a+t(b-a)}}^\downarrow(\bar t,1),
\]
But as $\of{\vp_{a,a+t(b-a)}}^\downarrow(\bar t,1)=\of{\vp_{a,b}}^\downarrow(\bar t,1)$ is true and by lower semicontinuity of $\vp_{a,b}$ the value $\vp(a+t(b-a))$ converges to $+\infty$ as $t$ converges to $t_0$, this implies that $\of{\vp_{a,b}}^\downarrow(\bar t,1)$ converges to $+\infty$ and eventually satisfies the desired inequality.

(b) Notice that $\vp_{a,b}(s)=\vp_{b,a}(1-s)$ and $\of{\vp_{a,b}}^\downarrow(s,-1)=\of{\vp_{b,a}}^\downarrow((1-s),1)$, hence the result is immediate from the above. \pend

\begin{corollary}\label{cor:scalar_Attainment_if_Deriv_leq_0}
Let $\vp \colon X \to \OLR$ be a radially l.s.c. function and $a \in \dom\vp$.
If $\vp^\downarrow(b, a-b) \leq 0$ for all $b \in X$, then either $\vp(a) = -\infty$ or $\vp$ is proper and $\dom \vp$ is star-shaped at $a$. In both cases, the infimum of $\vp$ is attained at $a$.
\end{corollary}

\proof
Theorem \ref{thm:Diewert} tells us that for all $b \in X$ there exists an $s \in (0,1]$ such that 
\[
\vp(a)\idif \vp(b)\leq \vp^\downarrow_{a,b}(s,-1).
\]
Using the definition of the lower Dini directional derivative one directly checks that
\[
\forall s \in \R \colon \vp^\downarrow_{a,b}(s,-1) = \vp^\downarrow(a + s(b-a), a-b)
\]
Taking $x = a + s\of{b-a}$ we obtain from $\vp^\downarrow(x, a-x) \leq 0$
\[
\forall s \in \R \colon \vp^\downarrow(a + s(b-a), a - \of{a + s(b-a)})  = \vp^\downarrow(a + s(b-a), -s(b-a)) \leq 0 
\]
Using the positive homogeneity of $\vp^\downarrow(x, \cdot)$ we get
\[
\forall s > 0 \colon \vp^\downarrow(a + s(b-a), a-b) =   \vp^\downarrow_{a,b}(s,-1) \leq 0.
\]
Hence $\vp(a)= -\infty$ or $-\infty < \vp(a) \leq \vp(b)$ for all $b\in X$. In the second case, $\vp$ is proper since $a \in \dom \vp$.

It is left to prove that $\dom \vp$ is star shaped at $a$.
Assume $b\in\dom \vp$ and $t\notin\dom\vp_{a,b}$ for some $t\in\of{0,1}$ and set
\[
r_0=\inf\cb{r\in\sqb{t,1}\st r\in\dom\vp_{a,b}}.
\]  
If $r_0\in\dom\vp_{a,b}$ then we are done, as in this case for $x=a+r_0(b-a)$  by lower semicontinuity of $\vp_{a,b}$ it holds $\vp^\downarrow(x,a-x)=+\infty$, a contradiction.
Hence assume $\vp_{a,b}(r_0)=+\infty$.
As $r_0<1$, we can chose a strictly decreasing sequence $\cb{r_n}_{n\in\N}\subseteq\dom\vp_{a,b}$ with $r_n\to r_0$ as $n$ converges to $+\infty$.
Applying Theorem \ref{thm:Diewert} to $a_n=a+r_{n+1}(b-a)$ and $b_n=a+r_n(b-a)$ for all $n\in\N$, then it exists a $0< t\leq 1$ such that for $r=r_{n+1}+t(r_n-r_{n+1})$ it holds
\[
\vp(a_n)\idif \vp(b_n)  \leq \of{\vp_{a_n,b_n}}^\downarrow(t,-1)=\of{\vp_{a,b}}^\downarrow(r,-1)
\]
Hence by assumption
\[
\vp(a+r_{n+1}(b-a))\idif \vp(a+r_{n}(b-a))  \leq \of{\vp_{a,b}}^\downarrow(r,-1)\leq 0,
\]
implying 
\[
\vp(a+r_{n+1}(b-a))  \leq \vp(a+r_{n}(b-a)).
\]
Especially, $\cb{\vp(a+r_n(b-a))}_{n\in\N}$ is a decreasing sequence in $\R\cup\cb{-\infty}$ as $\cb{r_n}_{n\in\N}\subseteq \dom \vp_{a,b}$ was assumed. By lower semicontinuity of $\vp_{a,b}$ it holds $\vp(a+r_0(b-a))\leq \liminf\limits_{n\to\infty}\vp(a+r_n(b-a))<+\infty$, a contradiction.
\pend

\medskip In the following definition, we extend some well-known notions to the case of extended real-valued functions, compare e.g. \cite{CambiniMartein2009, Crouzeix2005, GinchevIvanov2006, Ivanov01, ponstein1967}. Especially, we do not exclude the case $-\infty\in\vp\sqb{X}$ or  $\vp^\downarrow(b,a-b)=-\infty$.

\begin{definition}
\label{DefQuasiconvexEtc}
A function $\vp \colon X \to \OLR$ is said to be 

(a) quasiconvex if for all  $a,b\in X$ and all $t\in\of{0,1}$,
$\vp_{a,b}(t) \leq \max\cb{\vp(a),\vp(b)}$;

(b) semistrictly quasiconvex if for all  $a, b \in \dom\vp$ with $\vp(a) \neq \vp(b)$ and all $t \in \of{0,1}$,
$\vp_{a,b}(t) < \max\cb{\vp(a),\vp(b)}$;

(c) (lower Dini) pseudoconvex, if $\vp(a) < \vp(b)$ implies $\vp^\downarrow(b,a-b) < 0$;
\end{definition}

It is an easy task to prove that a convex function is semistrictly quasiconvex, quasiconvex and pseudoconvex.

Notice that semistrict quasiconvexity is defined with a strict inequality for all $a,b\in\dom \vp$  with $\vp(a) \neq \vp(b)$ while quasiconvexity only requires an inequality, but for all $a, b\in X$. 
The notions of a quasiconvex or semistrictly quasiconvex function are independent of each other as the following example shows.

\begin{example}\label{ex:semistr_quasiconv_NOT_quasiconv}
Let $\vp \colon \R \to \R$ be such that $\vp(0)=1$ and $\vp(x)=0$ for $x\neq 0$. Then $\vp$ is semistrictly quasiconvex, but not quasiconvex.
The function $\psi = -\vp$ is quasiconvex, but not semistrictly quasiconvex.
\end{example}

If $\vp \colon X \to \OLR$ is radially quasiconvex or semistrictly quasiconvex at $a \in \dom \vp$ then $\dom \vp$ is star-shaped at $a$.
The domain of an extended real-valued l.s.c. and pseudoconvex function is not necessarily star-shaped anywhere, therefore it does not have to be quasiconvex or semistriclty quasiconvex either. On the other hand, neither quasiconvexity, nor semistrict quasiconvexity implies pseudoconvexity, either.

\begin{example}
Let $\vp \colon \R \to\OLR$ be defined by $\vp(x)=0$ whenever $x \leq 0$ or $x \geq 1$ and $\vp(x)=+\infty$ otherwise. Then $\vp$ is  l.s.c. and pseudoconvex, but $\dom \vp$ is nowhere star-shaped, hence $\vp$ is neither quasiconvex, nor semistrictly quasiconvex.  
On the other hand, let $\psi\colon \R\to\OLR$ be defined as $\psi(x)=-x^2$, whenever $0\leq x$ and $\psi(x)=+\infty$, elsewhere. Then $\psi$ is both semistrict quasiconvex and quasiconvex, but $\psi^\downarrow(0,1)=0$ in contrast to $\psi(1)<\psi(0)=0$, hence $\psi$ is not pseudoconvex.
\end{example}

It is an easy task to prove that convexity of a function implies semistrict quasiconvexity, quasiconvexity and pseudoconvexity also for improper functions $\vp:X\to\OLR$. 

\begin{remark}\label{rem:qconvex_LevelSets_Convex}
The following equivalent characterizations of quasiconvexity are well known for proper functions, compare, for example, \cite[Proposition 3.2]{Crouzeix2005}. Without any problems, they can be extended to the general case of extended real-valued functions $\vp \colon X \to \OLR$.

(a1) The function $\vp \colon X \to \OLR$ is quasiconvex;

(a2) For all $r \in\OLR$ the lower level set $L_\vp(r)$ is convex;

(a3) For all $r \in \OLR$ the strict lower level set $L_\vp^<(r)=\cb{x\in X\st \vp(x)<r}$ is convex.

In particular, if $\vp$ is quasiconvex, then $\dom\vp$ and $L_\vp(-\infty)$ are convex sets.
\end{remark}

The following definition provides ``radial'' versions of the properties from Definition \ref{DefQuasiconvexEtc}.

\begin{definition}
\label{DefRadialQuasiconvexEtc}
A function $\vp \colon X \to \OLR$ is said to be radially quasiconvex (semistrictly quasiconvex, pseudoconvex) at $x_0 \in X$ if the function 
$\vp_{x_0, x} \colon \R \to \OLR$ is quasiconvex (semistrictly quasiconvex, pseudoconvex)  for all $x \in X$.
\end{definition}

\begin{proposition}\label{cor:Props_qconvex}
Let $\vp \colon X \to \OLR$ be a function. Then:

(a) The set $L^<_\vp(\vp(x))\cup\cb{x}$ is star-shaped at $x$ for all $x\in \dom \vp$ if, and only if, $\vp$ is semistrictly quasiconvex.

(b) If $\vp$ is semistrictly quasiconvex and l.s.c. then it is quasiconvex.

(c) A function $\vp$ is (semistrictly) quasiconvex if, and only if, it is radially (semistrictly) quasiconvex at every $x \in \dom \vp$.
\end{proposition}

\proof
(a) The function $\vp$ is semistrictly quasiconvex if, and only if, $\vp(y)<\vp(x)$ implies $\vp(y+t(x-y))<\vp(x)$ for all $t\in\of{0,1}$. This, in turn is equivalent to $L^<_\vp(\vp(x))\cup\cb{x}$ being star-shaped at $x$ for all $x \in X$.

(b) We only need to check the quasiconvexity inequality for $\vp(x) = \vp(y)$. Define $x_t = x + t(y-x)$ with $t \in \of{0,1}$ and assume $\vp(x_t) > \vp(x) = \vp(y)$. By semistrict quasiconvexity, $\vp(x_s) < \vp(x_t)$ for all $s \in [0, 1]\bs\cb{t}$. If $s \in (t, 1)$ and $\vp(x_s) \neq \vp(x)$ then again by semistrict quasiconvexity $\vp(x_t) < \max\cb{\vp(x),\vp(x_s)}$, a contradiction. The same can be done for $s \in \of{0,t}$, hence $\vp(x)=\vp(x_s)$ for all $s \in \of{0,1}\bs\cb{t}$ and $\vp(x_t) > \vp(x)$. This contradicts the lower semicontinuity of $\vp$. 

(c) Immediate. \pend

Especially, $\vp$ is radially semistrictly quasiconvex at $a\in\dom\vp$ if, and only if, for all $b \in \dom\vp$ and all $t \in \sqb{0,1}$ the set $L^<_{\vp_{a,b}}(\vp(a+t(b-a)))\cup\cb{t}$ is a convex interval. 

\begin{proposition}\label{cor:scalar_Deriv_leq_0_if_Attainment}
If $\vp(a) = \inf\vp\sqb{X}$ for some $a \in \dom\vp$ then $\vp$ is radially quasiconvex at $a$ if, and only if,
\begin{equation}\label{eq:qconvex_at_x}
\forall b \in X,\; \forall t \in \of{0,1} \colon \vp(a+t(b-a)) \leq \max\cb{\vp(a),\vp(b)}.
\end{equation}
In this case, $\vp^\downarrow(b,a-b)\leq 0$ holds true for all $b\in X$.
\end{proposition}

{\sc Proof.} If $\vp$ is radially quasiconvex at $a$, then \eqref{eq:qconvex_at_x}  is immediate. 

For the converse, let $\vp(a) = \inf\vp\sqb{X}$ and \eqref{eq:qconvex_at_x} be satisfied. Then $\vp(a)\leq\vp(x_t)$ is satisfied for all $b\in X$ and all $t\in\of{0,1}$ where $x_t = a + t(b-a)$. By \eqref{eq:qconvex_at_x}, $\vp(x_s)\leq \vp(x_t)$ for all $s\in\of{0,t}$. Now, take $s_1, s_2 \in [0,1]$ with $s_1 \neq s_2$, $\alpha \in (0, 1)$ and set $t = \max\cb{s_1, s_2}$, $s = \alpha s_1 + (1-\alpha)s_2$. Then $s \in (0, t)$ hence by the above
\[
\vp\of{x_s} \leq \vp\of{x_t} \leq \max\cb{\vp\of{x_{s_1}}, \vp\of{x_{s_2}}}
\]
which means that $\vp_{a, b}$ is quasiconvex since the remaining cases for $s_1, s_2, \alpha$ are trivial.

If the conditions of the first part are satisfied then $\vp(a)\leq \vp(b)$ hence, by \eqref{eq:qconvex_at_x}, $\vp(b+t(a-b)) \leq \vp(b)$ for all $t\in(0,1)$ which in turn implies $\vp(b+t(a-b))\idif \vp(b) \leq 0$  for all $t\in(0,1)$ whence
\[
\liminf\limits_{t\downarrow 0}\frac{1}{t}\of{\vp(b+t(a-b))\idif \vp(x)}\leq 0.
\]
This completes the proof. \pend

In general, Property \eqref{eq:qconvex_at_x} is weaker then radial quasiconvexity at $a$.

\begin{example}
Let $\vp\colon \R\to\OLR$ be given by $\vp(x)=\sup\cb{x^2,1-x^2}$. Then property \eqref{eq:qconvex_at_x} is satisfied at $a=-2$, but $\vp$ is not radially quasiconvex at $a$.
\end{example}

\begin{proposition}\label{lem:AttainmInfQuasiconvex} 
Let $\vp \colon X \to \OLR$ be radially l.s.c. at $a \in X$. Then $\vp$ is radially quasiconvex at $a$ if, and only if, for all $b \in X$ and all $r \in\OLR$ the set $\cb{t \in \sqb{0,1} \mid \vp_{a,b}(t) \leq r}$ is a closed convex subset of $\sqb{0,1}$ (a closed interval, possibly empty).

In this case,  the set $\cb{s \in \sqb{0,1} \mid \vp_{a,b}(s) = \inf\limits_{t \in \sqb{0,1}} \vp_{a,b}(t)}$ also is a closed convex subset of $\sqb{0,1}$ which is non-empty for each $b \in X$.
\end{proposition}

\proof
With Remark \ref{rem:qconvex_LevelSets_Convex} and the lower level set characterization of lower semi-continuity in view, the sublevel sets $L_{\vp_{a,b}}(r)$ are closed convex sets for all $b\in X$ and all $r\in\OLR$ if, and only if, the function $\vp_{a,b}$ is l.s.c. and quasiconvex for all $b \in X$. This proves the equivalence.

In this case, the set
\[
L_{\vp_{a,b}}\of{\inf\limits_{t \in \sqb{0,1}} \vp_{a,b}(t)}=\cb{s\in \sqb{0,1}\st \vp_{a,b}(s)=\inf\limits_{t \in \sqb{0,1}} \vp_{a,b}(t)}
\]
is closed and convex for each $b \in X$ which proves the second claim. This set is non-empty which is trivially the case if $-\infty$ is among the values of $\vp_{a,b}$, and which follows from the Weierstrass theorem since $\vp_{a,b}$ is lower semicontimuous on the compact set $[0,1]$.
\pend

\begin{proposition}\label{prop:strict_monotone}
Let $\vp \colon X \to \OLR$ be radially l.s.c. at $a \in \dom\vp$. Then $\vp$ is radially semistrictly quasiconvex at $a$ if, and only if, for all $b \in \dom\vp$ there exist $s_0 \leq t_0 \in \sqb{0,1}$ such that $\vp_{a, b}$ is strictly decreasing on $\sqb{0, s_0}$, strictly increasing on $\sqb{t_0,1}$ and constantly equal to $\inf\vp_{a,b}\sqb{0,1}$ on $\sqb{s_0, t_0}$.
\end{proposition}
\proof
Assume $\vp$ is radially semistrictly quasiconvex at $a \in \dom \vp$. Take $b \in \dom \vp$. By Proposition \ref{cor:Props_qconvex}, (b) $\vp$ is radially quasiconvex at $a$. Proposition \ref{lem:AttainmInfQuasiconvex} yields the existence of $s_0 \leq t_0 \in \sqb{0,1}$ such that $\sqb{s_0, t_0} = L_{\vp_{a, b}}\of{\inf\limits_{t\in\sqb{0,1}}\vp_{a, b}(t)}$. If $0 < t< s < s_0$ then $\vp_{a,b}(0) > \vp_{a,b}(t) > \vp_{a, b}(s)$ by semistrict quasiconvexity of $\vp_{a, b}$ and the fact that $s_0$ is a minimizer of $\vp_{a, b}$ on $[0,1]$. A similar argument proves that $\vp_{a,b}$ is strictly increasing on $\sqb{t_0, 1}$.

Conversely, let $0 \leq s < t \leq 1$ be such that $\vp_{a, b}(s) < \vp_{a, b}(t)$. Then $s, t \in [t_0, 1]$, hence $\vp_{a, b}\of{\alpha s + (1-\alpha)t} < \vp_{a, b}\of{t} = \max\cb{\vp_{a, b}\of{s}, \vp_{a, b}\of{t}}$ for all $\alpha \in (0, 1)$. If $0 \leq t < s \leq 1$ such that $\vp_{a, b}(s) < \vp_{a, b}(t)$ then $t, s \in [0, s_0]$ and a parallel argument works. Hence $\vp_{a, b}$ is semistrictly quasiconvex. \pend

\begin{proposition}\label{prop:pconv_lsc_then_quasiconv}
Let $\vp \colon X \to \OLR$  be radially pseudoconvex and radially l.s.c. at $a \in \dom\vp$ such that $\dom\vp$ is star-shaped at $a$. Then
$\vp$ is radially semistrictly quasiconvex at $a$. 
\end{proposition}

\proof Assume that for some $b\in\dom \vp$ the function $\vp_{a,b}$ is not semistrictly quasiconvex. Then there are $r, s, t \in \R$ such that
$0 \leq r < s < t \leq 1$, $\vp_{a,b}\of{r} \neq \vp_{a,b}(t)$ and
\[
\max\cb{\vp_{a,b}\of{r}, \vp_{a,b}(t)} \leq \vp_{a,b}(s).
\]
We assume $\vp_{a,b}\of{r} < \max\cb{\vp_{a,b}\of{r}, \vp_{a,b}\of{t}} = \vp_{a,b}\of{t}$. The other case can be dealt with by symmetric arguments.

Fix $\delta > 0$ such that $\vp_{a,b}\of{r} < \vp_{a,b}\of{t} - \delta$. Since $\vp_{a,b}$ is l.s.c. the set
\[
\cb{s' \in \R \mid \vp_{a,b}\of{s'} > \vp_{a,b}\of{t} - \delta}
\]
is open. Hence there is $\eps > 0$ such that $[s - \eps, s + \eps] \subseteq \of{r, t}$ and 
\[
\forall s' \in [s - \eps, s + \eps] \colon \vp_{a,b}\of{s'} > \vp_{a,b}\of{t}  - \delta. 
\]

Take $s' \in [s, s+\eps)$, $s'' \in (s', s+\eps]$ and assume $\vp_{a,b}\of{s''} < \vp_{a,b}\of{s'}$. By Theorem \ref{thm:Diewert} there exists an $\hat s \in (s' , s'']$ satisfying
\[
0 < \vp_{a,b}\of{s'} - \vp_{a,b}\of{s''} \leq \of{\vp_{a,b}}^\downarrow\of{\hat s, s' - s''}.
\]
Indeed, setting $a' = a+s'(b-a)$, $b' = a+s''(b-a)$ one obtains by Theorem \ref{thm:Diewert} an $\alpha \in (0, 1]$ satisfying $\vp\of{a'} - \vp\of{b'} \leq \of{\vp_{a,b}}^\downarrow\of{\alpha, -1}$. Defining $\hat s = s + \alpha(s''-s') \in (s', s'']$ and observing $\vp\of{a'} = \vp_{a,b}\of{s'}$, $\vp\of{b'} = \vp_{a,b}\of{s''}$ and $\of{\vp_{a,b}}^\downarrow\of{\alpha, -1} = \of{\vp_{a,b}}^\downarrow\of{\hat s, s' - s''}$ one obtains the above inequality. Using the positive homogeneity of the directional derivative we can multiply the inequality $0 < \of{\vp_{a,b}}^\downarrow\of{\hat s, s' - s''}$ by $\frac{r - \hat s}{s' - s''} > 0$ and obtain
$0 < \of{\vp_{a,b}}^\downarrow\of{\hat s, r - \hat s}$. The pseudoconvexity of $\vp_{a,b}$ yields $\vp_{a,b}\of{r} \geq \vp_{a,b}\of{\hat s}$ which contradicts the assumption $\vp_{a,b}\of{r} < \vp_{a,b}\of{t} - \delta < \vp_{a,b}\of{\hat s} - \delta$ (observe $\hat s \in [s, s+\eps]$). Hence $\vp_{a,b}\of{s''} \geq \vp_{a,b}\of{s'}$ whenever $s', s'' \in [s, s+\eps]$ and $s' < s''$. This implies
\[
\forall s' \in [s, s+\eps) \colon \of{\vp_{a,b}}^\downarrow\of{s', 1} \geq 0,
\]
and positive homogeneity of the directional derivative implies $\of{\vp_{a,b}}^\downarrow\of{s', t-s'} \geq 0$ and this by pseudoconvexity of $\vp_{a,b}$
\[
\vp_{a,b}\of{t} \geq \vp_{a,b}\of{s'} \geq \vp_{a,b}\of{s} \geq \vp_{a,b}\of{t}.
\]
This means $\vp_{a,b}\of{s'} = \vp_{a,b}\of{t}$ for all $s' \in [s, s+\eps)$. In turn, this implies that for $s' \in (s, s+\eps)$ we have $\of{\vp_{a,b}}^\downarrow\of{s', -1} \geq 0$, hence $\of{\vp_{a,b}}^\downarrow\of{s', r- s'} \geq 0$ and by pseudoconvexity $\vp_{a,b}\of{s'} \leq \vp_{a,b}\of{r}$. This contradicts the assumption $\vp_{a,b}\of{r} < \vp_{a,b}\of{t}$, hence (together with the symmetric case) the function $\vp_{a,b}$ is semistrictly quasiconvex for all $b \in \dom \vp$.\pend

By  Corollary \ref{cor:Props_qconvex}, a radially l.s.c. and radially semistrictly quasiconvex function $\vp:X\to\OLR$ especially is radially quasiconvex. Thus under the assumptions of Proposition \ref{prop:pconv_lsc_then_quasiconv} $\vp$ is also radially quasiconvex at $a$.

\begin{corollary}\label{cor:strict_monotone}
Let $\vp \colon X \to \OLR$  be radially pseudoconvex and radially l.s.c. at $a \in \dom\vp$ such that $\dom\vp$ is star-shaped at $a$. If $\vp^\downarrow(b,a-b)<0$ then $\vp^\downarrow(b_t, a-b_t)<0$  for all $t \geq 1$ where $b_t = a+t(b-a)$. If, additionally, $\vp(b)>-\infty$ then $\vp(b_t)>-\infty$ for all $t\geq 1$.
\end{corollary}
\proof The result is immediate if $\vp\of{b_t} = +\infty$ since in this case $\vp^\downarrow(b_t, a-b_t) = -\infty$ due to the properties of the inf-residuation $\idif$ on $\OLR$ and the definition of the directional derivative.

Assume $b_t \in \dom \vp$. Since $\vp^\downarrow(b,a-b)<0$ there exists an $s \in \of{0,1}$ such that either $\vp(a+s(b-a))<\vp(b)$ or $\vp(b)=\vp(a+s(b-a))=-\infty$. Hence, for $t > 0$ we either have $\vp(a+t(b-a))=\vp(b)=-\infty$ or, by Proposition \ref{prop:strict_monotone} applied to $\vp_{a, b_t}$, $\vp(a+t(b-a))>\vp(b)$. Note that, by Proposition \ref{prop:pconv_lsc_then_quasiconv}, $\vp$ is radially semistrictly quasiconvex at $a$. In both cases, $\vp^\downarrow(b_t,a-b_t)<0$ for $t \geq 0$ since in the first case we can apply that $\vp$ is radially pseudoconvex at $a$, and the second produces $\vp^\downarrow(b_t,a-b_t) = -\infty$ from $\of{-\infty} \idif \of{-\infty} = -\infty$.

Finally, if $\vp(b)>-\infty$, then, again by Proposition \ref{prop:strict_monotone} applied to $\vp_{a, b_t}$, $\vp(b_t)>\vp(b)$ for all $t>1$.
\pend

%%%New subsection
\subsection{Generalized convexity for set-valued functions}

In this section, we define (generalized) convexity notions for a set-valued function $f$, sometimes through the corresponding properties for the scalarizations $\vp_{f, z^*}$.

\begin{definition}
A function $f \colon X \to \G(Z,C)$ is called quasiconvex if
\begin{align}\label{eq:qconvex_sv}
\forall a, b \in X, \; \forall t \in \sqb{0,1} \colon f(a+t(b-a)) \supseteq f(a) \cap f(b).
\end{align}
\end{definition}

Formula \eqref{eq:qconvex_sv} is equivalent to $f(a+t(b-a)) \lel_C \sup\cb{f(a),f(b)}$ since the supremum in $\G(Z,C)$ is an intersection. Therefore, the definition of quasiconvexity for set-valued functions is a direct generalization of the scalar definition.

With respect to scalarizations we shall use the following concepts, compare \cite{BenoistBorwPop03, BenoistPopovici03, pjov6n1cre} and also the result presented in Theorem \ref{thm:Reverse_Attainment2} below.

\begin{definition}
\label{DefScalarizedRadialConvexities} A function $f \colon X \to \G(Z,C)$ is called 

(a) $\C$-l.s.c. if $\vp_{f, z^*}$ is l.s.c. for all $z^* \in \C$,

(b) radially $\C$-quasiconvex (semistrictly quasiconvex, pseudoconvex, l.s.c.) at $x_0 \in X$ if $\vp_{f,z^*} \colon X \to \OLR$ is radially quasiconvex (semistrictly quasiconvex, pseudoconvex, l.s.c.) at $x_0 \in X$ for all $z^*\in \C$.
\end{definition}

As in the scalar case, we introduce ``radial'' properties for set-valued functions as follows, compare  \cite{crespiJOTA}.

\begin{definition}
\label{DefRadialSetValued}
A function $f \colon X \to \G(Z,C)$ is called radially l.s.c. (radially quasiconvex) at $a \in X$ if the function $f_{a,b} \colon \R \to \G(Z,C)$ defined by
\begin{align*}%\label{eq:f_on_Intervall}
f_{a,b}(t)= \left\{\begin{array}{ccc}
			f(a+t(b-a)) & : & t \in \sqb{0,1}\\
						\emptyset & : &  \text{otherwise.}
			\end{array}
			\right.	
\end{align*}
is l.s.c. (quasiconvex).
\end{definition}

The equation
\[
\forall t \in \R \colon (\vp_{f,z^*})_{a,b}(t)=\vp_{f_{a,b},z^*}(t)
\]
is immediate.

Direct calculations prove that a set-valued function $f:X\to\G(Z,C)$ is convex, if and only if it is $\C$-convex, i.e. each scalarization $\vp_{f,z^*}:X\to\OLR$ with $z^*\in\C$ has a convex epigraph.
Moreover, a $\C$-quasiconvex function $f \colon X \to \G(Z,C)$ is quasiconvex, compare \cite[Theorem 2.1]{Kuroiwa96}. The following example shows that, in general, the second implication cannot be reversed. 
\begin{example}
Let $Z=\R^2$ and $C=\cl \co\of{ \cone\cb{(-1,1)^T,(1,1)^T}}$ and $f:\R\to Z$ be defined as
\[
f(x) = \left\{\begin{array}{lcl}
		(x,0)^T+C & : & x \in \pm2\mathbb{N} \\
		(2y+1,1)^T+C & : & 2y < x < 2(y+1), \, y\in\pm\mathbb{N}
	\end{array}
			\right.
\]
then $f$ is quasiconvex, while no scalarization with $z^*\in\Int C^-$ is quasiconvex.
\end{example}

\begin{remark} If $f \colon X \to \G(Z,C)$ is radially $\C$-pseudoconvex, radially $\C$-l.s.c. and $\dom f$ is star-shaped at $a$, then $f$ is radially $\C$-quasiconvex and radially $\C$-semistrictly quasiconvex at $a \in \dom f$. This follows from Proposition \ref{cor:Props_qconvex} (b) and Proposition \ref{prop:pconv_lsc_then_quasiconv}.
\end{remark}

\section{Characterization of infimizers}

According to the solution concept we introduced in Definition \ref{DefSetSolution}, we begin with the following definition.

\begin{definition}\cite{HamelSchrage12ArX}
Let $f \colon X \to \G(Z,C)$ and $M\subseteq X$ be non-empty. Then, the function $\hat f\of{\cdot; M} \colon X \to \G(Z,C)$ defined by
\begin{align*}
\hat f\of{x; M} = \inf f\sqb{M+x} = \inf_{m \in M} f\of{m + x} = \cl\co\bigcup_{m \in M}f\of{m + x}.
\end{align*}
is called the inf-translation of $f$ by $M$. The family of scalarizations of the inf-translation of $f$ by $M$ is given by
\[
\vp_{\hat f(\cdot; M),z^*}\of{x} = \inf_{z \in \hat f\of{x; M}}-z^*(z).
\]
\end{definition}

\begin{remark}
\label{RemInfTranslation}
The following relationships will be useful later on. We refer to \cite{HamelSchrage12ArX}.

(a) $\vp_{\hat f(\cdot; M),z^*}\of{x} = \hat\vp_{f,z^*}\of{x; M} = \inf_{m \in M}\vp_{f,z^*}\of{m + x}$.

(b) The infimum of $f\sqb{X}$ is attained in $M$, if and only if, it is attained in every $N \subseteq X$ with $M \subseteq N$; in particular, if $M$ is an infimizer then $\co M$ also is an infimizer.

(c) $\inf \hat f\of{\cdot; M}\sqb{X} = \inf\limits_{x \in X} \hat f\of{x; M} = \inf f\sqb{X}$.

(d) The infimum of $f\sqb{X}$ is attained in $M$, if and only if, $\hat f\of{0; M} = \inf \hat f\of{\cdot; M}\sqb{X}$.

(e) The infimum of $f\sqb{X}$ is attained in $M$, if and only if, $\vp_{f,z^*}\of{\cdot; M}$ attains its infimum at $0 \in X$ for all $z^*\in \C$,
\[
\inf f\sqb{M}=\inf f\sqb{X} \quad \Leftrightarrow \quad \forall z^*\in\C \colon \hat \vp_{f,z^*}\of{0; M} = \inf_{x \in X}\hat \vp_{f,z^*}\of{x; M}.
\]
This means that $0$ is a set $a$-minimizer of $f\of{\cdot; M}$ in the sense of \cite[Definition 3.2]{pjov6n1cre}, i.e. $\hat f\of{x; M} \subseteq \hat f\of{0; M}$ for all $x\in X$.
\end{remark}

\begin{proposition}\label{prop:UniformlyLscThenfM_lsc}
If $f \colon X \to \G(Z,C)$ is uniformly l.s.c. then $\hat f\of{\cdot; M}$ is $\C$-l.s.c. for all nonempty sets $M\subseteq X$.
\end{proposition}
\proof
If $f$ is uniformly l.s.c. then $\vp_{f,z^*}$ is uniformly l.s.c. for all $z^*\in \C$ as established in Proposition \ref{prop:Unif_lsc} (b). Replacing $x$ in \eqref{eq:uniform_lsc_scalar} by $m+x$ and taking the infimum over $m \in M$ on both sides yields that for every $\varepsilon>0$ there exists a $U\in \mathcal U_X(0)$ such that
\begin{align*}
\forall x\in X, \;\forall u\in U \colon \inf\limits_{m\in M}\vp_{f,z^*}(m+x)\leq  \inf\limits_{m\in M}\vp_{f,z^*}(m+x+u)+\varepsilon, 
\end{align*}
thus $\vp_{f\of{\cdot; M}, z^*} = \hat \vp_{f,z^*}\of{\cdot; M}$ is l.s.c.
\pend

The next result provides a sufficient condition for an infimizer in terms of the Dini directional derivative.

\begin{theorem}\label{thm:Attainment}
Let $f \colon X \to \G(Z,C)$ be uniformly l.s.c. and $\emptyset \neq M \subseteq \dom f$. If
\begin{align}\label{eq:strong_VI_2}
\forall x \in X \colon  \; 0 \in \bigcap\limits_{z^*\in \C}\hat f\of{\cdot; M}^\downarrow_{z^*}(x,-x)
\end{align}
then the infimum of $f$ over $X$ is attained in $M$. 
Moreover, $\hat f\of{0; M}=Z$ or $\dom \hat f\of{\cdot; M}$ is star-shaped at $0$.
\end{theorem}
\proof
By Proposition \ref{prop:UniformlyLscThenfM_lsc}, each scalarization $\vp_{f\of{\cdot; M}, z^*} = \hat \vp_{f,z^*}\of{\cdot; M}$ of $\hat f\of{\cdot; M}$ is (uniformly) l.s.c. Moreover, $0\in\dom \hat \vp_{f,z^*}\of{\cdot; M}$ since $M \subseteq \dom f$. From \eqref{EqDirDerScalarized} we conclude that \eqref{eq:strong_VI_2} is equivalent to $\vp_{f\of{\cdot; M},z^*}^\downarrow(x,-x)\leq 0$ for all $z^* \in \C$.

By Remark \ref{RemInfTranslation} (e), the infimum of $f\sqb{X}$ is attained in $M$ if, and only if, the infimum of $\vp_{f\of{\cdot; M},z^*}\sqb{X}$ is attained at $0$ for all $z^*\in \C$. 

Applying Corollary \ref{cor:scalar_Attainment_if_Deriv_leq_0} we obtain the results.
\pend

\begin{remark}
The condition \eqref{eq:strong_VI_2} in Theorem \ref{thm:Attainment} can be replaced by
\[
\hat f\of{0; M} = \hat f\of{0; \co M} \quad \text{and} \quad \forall x \in X \colon 0\in \bigcap\limits_{z^*\in \C}\hat f\of{\cdot; \co M}^\downarrow_{z^*}(x,-x).
\]
In this case, $\hat f\of{0; M} = \hat f\of{0; \co M} = Z$, or $\dom \hat f\of{\cdot; \co M}$ is star-shaped at $0$.
\end{remark}

\begin{remark} If $M = \cb{x_0}$ for $x_0 \in X$ then \eqref{eq:strong_VI_2} is equivalent to
\[
\forall x \in X \colon 0\in \bigcap\limits_{z^*\in \C}f^\downarrow_{z^*}\of{x, x_0 - x}.
\]
Indeed, this follows from $\hat f\of{x; \cb{x_0}} = f\of{x_0 + x}$ and, especially, $\hat f\of{0; \cb{x_0}} = f\of{x_0}$. Thus, if the infimizer is a singleton then the complicated looking condition \eqref{eq:strong_VI_2} boils down to a more familiar form. Although it is in general very unlikely that the infimum of a $\G(Z,C)$-valued function is attained in a single point, this is the case for the inf-translation of $f$ by an infimizer (set) $M$. The reduction of infimizer sets to singletons was the main motivation for the introduction of the inf-translation in \cite{HamelSchrage12ArX}.
\end{remark}

\begin{lemma}\label{lem:quasiconvex_extension}
Let $f \colon X \to \G(Z,C)$, $\emptyset \neq M \subseteq \dom f$ and $z^* \in \C$. Assume that \eqref{eq:qconvex_at_x} is satisfied for $\vp_{f, z^*}$ whenever $x_0 \in \co M$. Then \eqref{eq:qconvex_at_x} with $x_0 = 0$ is satisfied for $\hat \vp_{f, z^*}\of{\cdot; \co M}$.
\end{lemma}
\proof
Assume there are $x \in X$ and $t \in \of{0,1}$ such that $\hat \vp_{f, z^*}\of{\cdot; \co M}$ does not satisfy \eqref{eq:qconvex_at_x} with $x_0 = 0$, i.e.
\[
\hat \vp_{f,z^*}\of{tx; \co M} > \max\cb{\hat \vp_{f,z^*}\of{0; \co M}, \hat \vp_{f,z^*}\of{x; \co M}}.
\]
Since
\[
\hat \vp_{f,z^*}\of{y; \co M} = \inf_{m \in \co M} \vp_{f, z^*}\of{m + y}
\]
there are $m_1, m_2 \in \co M$ such that
\[
\forall m \in \co M \colon \vp_{f,z^*}\of{m + tx} > \max\cb{\vp_{f,z^*}\of{m_1}, \vp_{f,z^*}\of{m_2 + x}}.
\]
Taking $m  = m_1 + t\of{m_2 - m_1} \in \co M$ we obtain
\[
\vp_{f,z^*}\of{m_1 + t\of{x + m_2 - m_1}} > \max\cb{\vp_{f,z^*}\of{m_1}, \vp_{f,z^*}\of{m_2 + x}}
\]
which contradicts the assumption that $\vp_{f,z^*}$ satisfies \eqref{eq:qconvex_at_x} at any $x_0 \in \co M$ (choose $x_0 = m_1$ and replace $x$ in \eqref{eq:qconvex_at_x} in by $x + m_2$ with $x$ from above).
\pend

Combining the previous results we obtain the following necessary condition for infimizers.

\begin{theorem}\label{thm:Reverse_Attainment2}
Let $f \colon X \to \G(Z,C)$ and $\emptyset \neq M \subseteq \dom f$ be such that for each $z^*\in \C$ the scalarization $\vp_{f, z^*}$ of $f$ satisfies \eqref{eq:qconvex_at_x} whenever $x_0 \in M$. If the infimum of $f$ over $X$ is attained in $M$ then $\hat f\of{0; M} = \hat f\of{0; \co M}$, $\hat f\of{\cdot; \co M}$ is radially quasiconvex at $0$ and
\begin{align}\label{eq:Strong_VI2}
\forall x \in X \colon 0 \in \bigcap\limits_{z^*\in \C}\hat f\of{\cdot; \co M}^\downarrow_{z^*}(x,-x).
\end{align}
\end{theorem}

{\sc Proof.} Using Remark \ref{RemInfTranslation} (d), (b) we obtain $\hat f\of{0; M} = \hat f\of{0; \co M}$. Remark \ref{RemInfTranslation} (e), Lemma \ref{lem:quasiconvex_extension} and Proposition \ref{cor:scalar_Deriv_leq_0_if_Attainment} yield that $\hat \vp_{f, z^*}\of{\cdot; \co M}$ is radially quasiconvex at $0$, hence $\hat f\of{\cdot; \co M}$ is radially quasiconvex (see discussion after Definition \ref{DefRadialSetValued}). The derivative conditions now follows from Proposition \ref{cor:scalar_Deriv_leq_0_if_Attainment} and \eqref{EqDirDerScalarized}. \pend

\begin{remark}
Notice that radial quasiconvexity of each scalarization of $f$ at each $m\in \co M$  and $\inf f\sqb{M} = \inf f\sqb{X}$ together are sufficient %, but not neccessary 
conditions for the assumptions of Theorem \ref{thm:Reverse_Attainment2} to be satisfied.
\end{remark}

If $Z$ is a Banach space and $C+(-C)=Z$, i.e. $C$ generates $Z$, the function $f \colon X \to \G(Z,C)$ is quasiconvex and the infimum of $f\sqb{X}$ is attained in $\co M\subseteq \dom f$, then 
\begin{align}\label{eq:Strong_VI_extd}
0\in\bigcap\limits_{z^* \in \extd C^-}\hat f\of{\cdot; \co M}^\downarrow_{z^*}(x,-x).
\end{align}
However, \eqref{eq:Strong_VI_extd} can hold without $f\of{0; \co M}$ being anywhere near the infimum of $f\sqb{X}$.
Therefore, the sufficient property given in Theorem  \ref{thm:Attainment} therefore is notably stronger. 
If $f \colon X \to \G(Z,C)$ is uniformly l.s.c. and \eqref{eq:Strong_VI_extd} is satisfied, then $0$ is a set $A$-minimizer of $f\of{\cdot; \co M}$ in the sense of \cite[Definition 3.4]{pjov6n1cre}, i.e. for all $x\in X$ and all $z^*\in\extd C^-$ it holds $\vp_{\hat f\of{\cdot; \co M},z^*}(0)\leq \vp_{\hat f\of{\cdot; \co M},z^*}(x)$.

\begin{example}
Let $Z=\R^2$, $C=\R^2_+$ and $f=\R\to\G(Z,C)$, $\dom f=\R_+$ and $f(x)=\sqb{(0,x)^T,(x,0)^T}+C$, whenever $x\in\dom f$. 
Thus $Z$ is a Banach space, $C$ generates $Z$, $f$ is uniformly l.s.c. and convex and thus especially $\C$--quasiconvex, the infimum of $f\sqb{X}$ and $\vp_{f,z^*}\sqb{X}$ is attained in $0$, whenever $z^*\in\Int C^-$.
The extreme directions of $C^-$ are the elements of $\cone\cb{(-1,0)^T,(0,-1)^T}$ and whenenver $z^*\in\extd C^-$, then $\vp_{f,z^*}\sqb{\dom f}=\cb{0}$.  
Thus, for any choice of $\emptyset\neq M\subseteq \dom f$, \eqref{eq:Strong_VI_extd} is satisfied whereas \eqref{eq:strong_VI_2} holds true, iff $0\in M$.
\end{example}

%%%New section
\section{Characterization of minimizers}\label{sec:NondomElem}

In this section, we shall give sufficient conditions for a point $x_0 \in \dom f$ to produce a minimal value of $f \colon X \to \G(Z,C)$, i.e. $f(x_0)$ satisfies
\[
x \in X, \; f\of{x} \supseteq f\of{x_0} \quad \Rightarrow \quad f\of{x} = f\of{x_0}
\]
(compare (b) of Definition \ref{DefSetSolution}).

For a fixed $x_0 \in \dom f$, we define the set
\[
A(f,x_0) = \cb{x \in X \mid f(x) \not\subseteq f(x_0)}
\]
which, of course, always is a subset of $\dom f$.  Note 
\begin{align*}
f(x_0) \in \Min f\sqb{X} \quad \Leftrightarrow \quad  \sqb{x \in X, \; f\of{x} \supseteq f\of{x_0} \; \Rightarrow \; x \not\in A(f,x_0)}.
\end{align*}
 If $x \in A(f,x_0)$ we can separate a point $\bar z \in f\of{x} \bs f\of{x_0}$ from $f\of{x_0}$ since the images of $f$ are closed convex sets. Thus, there are $z_0^*\in \C$, $r_0 \in \R$ such that
\[
\forall z \in f\of{x_0} \colon z^*_0\of{\bar z} < r_0 \leq z^*_0\of{z}.
\]
Therefore, $x \in A(f,x_0)$ if, and only if,
\begin{align*}
\exists z_0^*\in \C \colon \vp_{f,z_0^*}(x) < \vp_{f,z_0^*}(x_0).
\end{align*}
Hence, if $A(f,x_0) \neq \emptyset$ then there is $z^*_0 \in \C$ such that $\vp_{f, z_0^*}(x_0) \in \R$ and $A(f,x_0) \subseteq \dom f$.

This discussion can be used to verify the following result.

\begin{proposition}
\label{PropAStar-Shaped} Let $f \colon X \to \G(Z,C)$ be radially $\C$-semistrictly quasiconvex at $x_0 \in \dom f$. Then $A\of{f, x_0}\cup\cb{x_0}$ is star-shaped at $x_0$.
\end{proposition}

\proof Assume there are $x \in A\of{f, x_0}$ and $t \in (0, 1)$ such that $x_0 + t(x-x_0) \not\in A\of{f, x_0}$. Then $f\of{x_0 + t(x-x_0)} \subseteq f\of{x_0}$, hence
\[
\forall z^* \in \C \colon \vp_{f, z^*}\of{x_0} \leq \vp_{f, z_0^*}\of{x_0 + t(x-x_0)}.
\]
On the other hand, the above separation argument shows $\vp_{f,z_0^*}(x) < \vp_{f,z_0^*}(x_0)$, hence by semistrict quasiconvexity of $\vp_{f, z_0^*}$
\[
\vp_{f, z_0^*}\of{x_0 + t(x-x_0)} < \max\cb{\vp_{f, z_0^*}\of{x_0}, \vp_{f, z_0^*}\of{x}} = \vp_{f, z_0^*}\of{x_0},
\]
a contradiction. \pend

We will prove that if a certain variational inequality of Minty type is satisfied for all $x \in A(f,x_0)$, then $f(x_0)$ is a minimal element in $f\sqb{X}$.

\begin{theorem}\label{thm:VI_Principle}
Let $f \colon X \to \G(Z,C)$ be radially $\C$-l.s.c. and radially $\C$-semistrictly quasiconvex at $x_0 \in \dom f$. If there is a non-empty finite set $M^* \subseteq \C$ such that
\begin{align}\label{eq:VI}
\forall x \in A(f,x_0), \; \exists z^*\in M^* \colon 0 \in \Int f^\downarrow_{z^*}(x,x_0-x) \, \wedge \, \vp_{f,z^*}(x)\neq-\infty
\end{align}
then $f(x_0) \in \Min f\sqb{X}$.
\end{theorem}

\proof
If $A(f,x_0) = \emptyset$ then $f(x) \subseteq f(x_0)$ for all $x\in X$, hence $f(x_0)=\inf f\sqb{X}$ and especially $f(x_0)\in\Min f\sqb{X}$. 

Assume $A(f,x_0) \neq \emptyset$ and fix $x \in A(f, x_0)$, $z^* \in \C$. Since $f$ is radially $\C$-semistrictly quasiconvex at $x_0$, Proposition \ref{PropAStar-Shaped} ensures that $A(f, x_0)\cup\cb{x_0}$ is star-shaped at $x_0$.

Since $A(f, x_0) \subseteq \dom f$ Proposition \ref{prop:strict_monotone} yields that for all $m^*\in M^*$ the value
\[
t(m^*) = \max\cb{t \in \sqb{0,1} \mid \of{\vp_{f,m^*}}_{x_0,x} \; \text{is non-increasing on} \; \sqb{0,t}}
\]
is well defined, $\vp_{f, m^*}(x_0 + t(m^*)(x - x_0)) = \inf \vp_{f, m^*}\sqb{x_0,x}$ and for all $t \in \sqb{0, t(m^*)}$ either $\vp_{f,m^*}(x_t) = -\infty$ or $\vp^\downarrow_{f,m^*}(x_t, x_0 - x_t)\geq 0$ where $x_t = x_0 + t(x-x_0)$. Since $M^*$ is finite, there exists $m^*_0 \in M^*$ such that
\[
t(m^*_0) = t_0 = \min\cb{t(m^*) \mid m^* \in M^*}.
\]
For all $m^*\in M^*$, either $\vp_{f,m^*}(x_{t_0}) = -\infty$, or $\vp^\downarrow_{f,m^*}(x_{t_0}, x_0 - x_{t_0}) \geq 0$. Since $A(f, x_0) \cup \cb{x_0}$ is star-shaped at $x_0$ the whole line segment between $x$ and $x_0$ belongs to $A\of{f, x_0}$ and therefore \eqref{eq:VI} is in force. The scalarization of the Dini derivative \eqref{EqDirDerScalarized} tells us that $0 \in \Int f^\downarrow_{m^*}(x,x_0-x)$ is equivalent to $\vp_{f,m^*}^\downarrow(x,x_0 - x) < 0$. Hence  \eqref{eq:VI} implies $t_0 = 0$.

Thus by Proposition \ref{prop:strict_monotone} $\of{\vp_{f,m^*_0}}_{x_0,x}$ is strictly increasing on $\sqb{0,1}$ which implies
\[
\vp_{f,m^*_0}(x_0)<\vp_{f,m^*_0}(x).
\]
According to Remark \ref{RemSolSetValuedScalarized}, (b) this verifies $f(x_0)\nsubseteq f(x)$ for all $x \in A(f, x_0)$. Finally, if $x\notin A(f, x_0)$ then $f(x) \subseteq f(x_0)$, hence $f(x_0)$ is minimal in $f\sqb{X}$.
\pend

\medskip We have proven that under the assumptions of Theorem \ref{thm:VI_Principle}, 
to any ray $\sqb{x,x_0}\subseteq \dom f$ with $\sqb{x,x_0}\cap A\neq \emptyset$
there exists a single element $m^*\in M^*$ such that Property \eqref{eq:VI} is satisfied for all $x_t=x_0+t(x-x_0)$ with $0<t\leq 1$.

\begin{remark}
A sufficient condition for radial $\C$-semistrict quasiconvexity  of $f$ in Theorem \ref{thm:VI_Principle} reads as follows. Let $f \colon X \to \G(Z,C)$ be radially $\C$-l.s.c and radially $\C$-pseudoconvex at $x_0 \in \dom f$ such that $\dom f$ is star-shaped at $x_0$. Then $f$ is radially $\C$-semistrictly quasiconvex at $x_0$. Indeed, in this case $\of{\vp_{f,z^*}}_{x_0, x}$ is semistrictly quasiconvex by Proposition \ref{prop:pconv_lsc_then_quasiconv} since, by Definition \ref{DefScalarizedRadialConvexities}, it is l.s.c., pseudoconvex and $\dom\of{\vp_{f,z^*}}_{x_0, x} \subseteq [0,1]$ is an intervall (including $0$) because $\dom f = \dom \vp_{f,z^*}$ is star-shaped at $x_0$ by assumption. Hence $f$ is radially $\C$-semistrictly quasiconvex.
\end{remark}

The following example shows that the assumption $M^*\subseteq \C$ be finite cannot be relaxed.

\begin{example}
Define $z^*_i=-\frac{1}{i+1}(1,i)^T\in (R^2_+)^-\bs\cb{0}$ for all $i \in \N = \cb{0,1, 2, \ldots}$. Let $f \colon \R \to \G(\R^2,\R^2_+)$ be such that $\dom f=\sqb{0,1}$ and $f(x)=\bigcap\limits_{i \in \N}L_{z^*_i}(-\vp_{z^*_i}(x))$ for all $x\in\sqb{0,1}$ where
\[
\vp_{z^*_i}(x)=\left\{\begin{array}{lcl}
-(i+1)\min\cb{1-x,i x}&:&\text{ if $x\in\sqb{0,1}$ and $i\in\N$;}\\
+\infty&: &\text{ elsewhere.}
	\end{array}
			\right.
\]
Since $\vp_{z^*_i}$ is a convex l.s.c function for all $i\in\N$, $\gr f$ is closed and convex, hence $f$ is l.s.c. and convex, and it is easy to see that $f(0)=f(1)=\R^2_+$. Defining $z_i(x)\in\R^2$ by 
\begin{align*}
\forall i \in \N\bs\cb{0} \colon 
\cb{z_i(x)}=\cb{z\in Z\colon z^*_{i-1}(z)=\vp_{z^*_{i-1}}(x)}\cap \cb{z\in Z\colon z^*_{i}(z)=\vp_{z^*_{i}}(x)}
\end{align*}
then $f(x)=\co\cb{z_i(x)\st i \in \N\bs\cb{0}} + C$ is true for all $x\in\of{0,1}$. This implies that $\vp_{f,z^*_i}(x)=\vp_{z^*_i}(x)$ is true for all $x\in\sqb{0,1}$ and all $i\in \N$ and therefore $f(x)\supsetneq f(0)$ is satisfied for all $x\in\of{0,1}$ and $f(0)\notin\Min f\sqb{\R}$.

On the other hand, for any given $x\in\of{0,1}$, it exists an $i \in \N\bs\cb{0}$ such that $x\in\of{\frac{1}{i+1},1}$, hence
$\vp^\downarrow_{f,z^*_i}(x,0-x)=-(i+1)<0$ and $-i\leq \vp_{f,z^*_i}(x)\neq-\infty$.
Hence the assumptions of Theorem \ref{thm:VI_Principle} are satisfied for $x_0=0$, replacing the finite set $M^*$ by $\C$, while $f(0) \notin\Min f\sqb{\R}$.
\end{example}

\begin{remark}
Recall that an element $z^* \in \C$ is an extreme direction of $C^-$ if for all $z^*_1, z^*_2 \in C^-$, $z^* = z^*_1 + z^*_2$ implies $z^*_1 = tz^*$ and $z^*_2 = sz^*$ for some nonnegative $t, s \in \R$, and the set of all extreme directions of $C^-$ is denoted by $\extd C^-$.

Let $M^*\subseteq\extd C^-$ be a nonempty, finite set with $\extd C^-=\cone M^*$. If the assumptions of Theorem \ref{thm:VI_Principle} are satisfied for this set $M^*$, then $f(x_0)\in\Min f\sqb{X}$. However we do not restrict the choice of the set $M^*$ to any specific subset of $\C$ thus the result of Theorem \ref{thm:VI_Principle} is true in a more general case, too.
Notice that $\extd C^-\neq\emptyset$ is rather restrictive, as for example it excludes such cases where $C^-$ contains linear subspaces of $Z^*$. 
\end{remark}

\begin{corollary}\label{cor:VI_Principle}
Let $f \colon X \to \G(Z,C)$ and $M \subseteq \dom f$ be such that $f$ is radially $\C$-l.s.c. and radially $\C$-semistrictly quasiconvex at every $u \in M$. Further, let $M^*\subseteq\C$ be a nonempty finite set. If for all $u \in M$ and for all $x \in X$ either $f(x) \subseteq f(u)$ or
\begin{align*}
\exists z^*\in M^* \colon 0 \in \Int f^\downarrow_{z^*}(x, u-x) \wedge \vp_{f,z^*}(x) \neq -\infty
\end{align*}
then $f\sqb{M} \subseteq \Min f\sqb{X}$.
\end{corollary}
\proof
The assumptions guarantee that for each $x \in X$ either $x \not\in A\of{f, u}$ for all $u \in M$ or $x \in A\of{f, u}$ and the condition in \eqref{eq:VI} is satisfied. Theorem \ref{thm:VI_Principle} produces the result.
\pend

The next result provides a necessary condition for a minimizer in terms of the Dini directional derivative.

\begin{theorem}\label{thm:Reverse_Nondom}
Let $f \colon X \to \G(Z, C)$ and $x_0 \in \dom f$ be such that $f(x_0)\in\Min(f\sqb{X})$. Assume $\vp_{f, z^*}$ satisfies \eqref{eq:qconvex_at_x} for all $z^*\in \C$. Then

(a) for all $x \in X$ there exists $z^* \in \C$ such that
\begin{equation}\label{eq:Reverse_Nondom_a}
0\in f^\downarrow_{z^*}(x, x_0-x),
\end{equation}

(b) if, additionally, $f$ is radially $\C$-l.s.c. and radially $\C$-pseudoconvex at $x_0$, then for all $x\in X$ either $f(x_0)=f(x)$, or there exists $z^*\in\C$ such that
\begin{equation}\label{eq:Reverse_Nondom_b}
0\in \Int f^\downarrow_{z^*}(x,x_0-x) \wedge \vp_{f,z^*}(x)\neq-\infty.
\end{equation}
\end{theorem}
\proof (a) According to Remark \ref{RemSolSetValuedScalarized}(b), for each $x \in X$ there exists $z^*_0 \in \C$ such that 
\[
\max\cb{\vp_{f, z^*_0}(x_0), \vp_{f, z^*_0}(x)} = \vp_{f, z^*_0}(x).
\]
Condition \eqref{eq:qconvex_at_x} now implies
\[
\forall t\in\sqb{0,1} \colon \vp_{f,z^*_0}(x+t(x_0-x)) \leq \vp_{f,z^*_0}(x)
\]
which in turn yields $\vp_{f,z^*_0}^\downarrow(x,x_0-x) \leq 0$. Equation \eqref{EqDirDerScalarized} produces the result.

(b) Under the additional assumption, $f$ is radially $\C$-semistrictly quasiconvex by Proposition \ref{prop:pconv_lsc_then_quasiconv}, thus for all $x\in X$
either $f(x+t(x_0-x)) = f(x)$ for all $t\in \sqb{0,1}$, or there exists $s \in \of{0,1}$ and $z^*_0 \in \C$ such that $\vp_{f,z^*_0}$ is strictly increasing on $\sqb{x+s(x_0-x), x}\cap\dom f$ by Proposition \ref{prop:strict_monotone}. Hence  $\vp_{f,z^*_0}^\downarrow(x,x_0-x) < 0$  by pseudoconvexity of $\vp_{f, z^*_0}$.
\pend

Stating the assumptions of Theorem \ref{thm:Reverse_Nondom} for all elements of a set $M\subseteq \dom f$, the following corollary is straightforward.

\begin{corollary}\label{cor:Reverse_Nondom}
Let $f \colon X \to \G(Z, C)$ and $\emptyset \neq M \subseteq \dom f$ be such that $f\sqb{M} \subseteq \Min(f\sqb{X})$. Assume that for all $z^* \in \C$ the function $\vp_{f, z^*}$ satisfies \eqref{eq:qconvex_at_x} with $x_0$ replaced by an arbitrary $u \in M$. Then

(a) for all $u \in M$ and all $x \in X$ there exists $z^* \in \C$ such that
\begin{equation}\label{eq:Reverse_Nondom_a}
0\in f^\downarrow_{z^*}(x, u - x),
\end{equation}

(b) if, additionally, $f$ is radially $\C$-l.s.c. and radially $\C$-pseudoconvex at $u$ for all $u \in M$, then for $u \in M$ and $x \in X$ either $f(u) = f(x)$, or 
\begin{equation}\label{eq:Reverse_Nondom_b}
\exists z^* \in \C \colon 0 \in \Int f^\downarrow_{z^*}(x, u - x) \; \wedge \; \vp_{f,z^*}(x) \neq -\infty.
\end{equation}
\end{corollary}

%%%New section
\section{Conclusions}

The combination of Theorem \ref{thm:Attainment} and Corollary \ref{cor:VI_Principle} produces the following sufficient condition for solutions of our basic set-valued optimization problem, i.e. of
\[
\tag{P} \text{minimize} \quad f \quad \text{subject to} \quad x \in X.
\]
Note that a set $M \subseteq X$ is a solution of \eqref{eq:P} if the infimum of $f\sqb{X}$ is attained in $M$ and $f\sqb{M}\subseteq \Min(f\sqb{X})$.

\begin{theorem}\label{thm:Concl_1}
Let $M^* \subseteq \C$ be a finite set, $f \colon X \to \G(Z,C)$ be a uniformly l.s.c. function and $\emptyset \neq M \subseteq \dom f$ such that $f$ is radially $\C$-semistrictly quasiconvex at $u$ for all $u \in M$. Moreover, let
\begin{align*}
\forall x \in X \colon 0 \in \bigcap\limits_{z^*\in \C}\hat f(\cdot,M)^\downarrow_{z^*}(x,-x)
\end{align*}
be satisfied and for $u \in M$ and $x \in X$ either $f(x)\subseteq f(u)$ or
\begin{align*}
\exists z^*\in M^* \colon 0\in\Int f^\downarrow_{z^*}(x, u - x) \wedge \vp_{f,z^*}(x)\neq-\infty.
\end{align*}
Then, $M$ is a solution of \eqref{eq:P}.
\end{theorem}
\proof A uniformly l.s.c. functions $f \colon X \to \G(Z,C)$ is (uniformly) $\C$-l.s.c. by Proposition \ref{prop:Unif_lsc}. The result follows from Theorem \ref{thm:Attainment} and Corollary \ref{cor:VI_Principle}.
\pend

Likewise, the combination of Theorem \ref{thm:Reverse_Attainment2} and Corollary \ref{cor:Reverse_Nondom} produces a necessary optimality condition for solutions of \eqref{eq:P}.

\begin{theorem}\label{thm:Concl_2}
Let $f \colon X \to \G(Z,C)$ and $\emptyset \neq M \subseteq \dom f$ be such that $M$ is a solution of \eqref{eq:P}. Assume that for all $z^* \in \C$ the function $\vp_{f, z^*}$ satisfies \eqref{eq:qconvex_at_x} with $x_0$ replaced by an arbitrary $u \in M$. 
Then
\begin{align*}
\forall x \in X \colon 0 \in \bigcap\limits_{z^*\in \C}\hat f(\cdot, \co M)^\downarrow_{z^*}(x, -x)
\end{align*}
is satisfied. Moreover, $\hat f(0, M) = \hat f(0, \co M)$ and 
\begin{equation*}
\forall u \in M, \; \exists z^* \in \C \colon 0\in f^\downarrow_{z^*}(x, u - x).
\end{equation*}

If, additionally, $f$ is radially $\C$-l.s.c. and radially $\C$-pseudoconvex at $u$ for all $u \in M$, then for $u \in M$ and $x \in X$ either $f(u) = f(x)$, or 
\begin{equation}\label{eq:Reverse_Nondom_b}
\exists z^* \in \C \colon 0 \in \Int f^\downarrow_{z^*}(x, u - x) \; \wedge \; \vp_{f,z^*}(x) \neq -\infty.
\end{equation}
\end{theorem}

\proof This directly follows from Theorem \ref{thm:Reverse_Attainment2} and Corollary \ref{cor:Reverse_Nondom}. \pend

\end{document}